\documentclass[a4paper,11pt]{amsart}
\usepackage{amsmath}%
\usepackage{amsfonts}%
\usepackage{amssymb}%
\usepackage{graphicx}%

\usepackage{comment}%    for commenting out when revising
\usepackage{tikz}%			to make pictures  
\usepackage{caption}%    for captioning without number
\usepackage{mathtools} % for \coloneqq
\usepackage{scalerel,stackengine} %for reallywidehat

\theoremstyle{plain}
\newtheorem{theorem}{Theorem}

\theoremstyle{definition}
\newtheorem{definition}[theorem]{Definition}

\theoremstyle{remark}

\theoremstyle{remark}
\newtheorem{remark}[theorem]{Remark}
 
\theoremstyle{remark}
\newtheorem{example}[theorem]{Example}

\newcommand{\mathsc}[1]{\textnormal{\textsc{#1}}} % so textsc behaves in formulas when in bold; reverts to nonbold
\newcommand{\Xb}{X_{\mathsc{Ban}}}

\newcommand{\CC}{\mathbf{C} }
\newcommand{\DD}{\mathbf{D} }

\newcommand{\PP}{\mathbf{P} }
\newcommand{\RR}{\mathbf{R} }

\newcommand{\ZZ}{\mathbf{Z}\,}
\newcommand{\CCC}{\mathcal{C} }

\newcommand{\FF}{\mathcal{F} }

\newcommand{\OO}{\mathcal{O} }

\newcommand{\Coh}{\operatorname{Coh}}

\newcommand{\HH}{\operatorname{H}}

\newcommand{\Pic}{\operatorname{Pic}}

\newcommand{\Bl}{\operatorname{Bl}}

\newcommand{\Jphi}[2]{\phi_{-2,1}(#1,#2)}
\newcommand{\Jphiqp}{\phi_{-2,1}(q,p)}

\newcommand{\JphiQ}[1]{\phi_{\scriptscriptstyle Q}(#1)}
\newcommand{\thetaQ}[1]{\theta_{\scriptscriptstyle Q}(#1)}
\newcommand{\etaQ}{\eta_{\scriptscriptstyle Q}}

\newcommand{\tor}{\CC^*\times\CC^*}

\newcommand{\Fb}{{F_{\mathsc{ban}}}}
\newcommand{\Xm}{X_{\mathsc{mb}}}
\newcommand{\Xd}{X_{\mathsc{22}}}

\newcommand{\Fd}{F_{\mathsc{22}}}

\newcommand{\Xq}{X_{\mathsc{15}}}

\newcommand{\Fm}{{F_{\mathsc{mb}}}}
\newcommand{\Fvw}[2]{F_{\mathsc{mb}}^{#1#2}}

\newcommand{\Fbhat}{{\widehat{F}_{\mathsc{ban}}}}
\newcommand{\Fmhat}{{\widehat{F}_{\mathsc{mb}}}}

\newcommand{\Lm}{\widehat{F}_{\mathsc{mb}}}
\newcommand{\Lvw}[2]{{\Lm}^{#1#2}}

\newcommand{\Ld}{\widehat{F}_{\mathsc{22}}}

\newcommand{\UFb}{U({F_{\mathsc{ban}}})}
\newcommand{\UFm}{U({F_{\mathsc{mb}}})}

\newcommand{\UFmhat}{U({\widehat{F}_{\mathsc{mb}}})}

\newcommand{\MF}{M_{\beta}^{\Fmhat}}

\newcommand{\naivel}{\widetilde{n}^0_{\beta}(\Lm)}

\makeatletter
\newcommand*\bigcdot{\mathpalette\bigcdot@{.5}}
\newcommand*\bigcdot@[2]{\mathbin{\vcenter{\hbox{\scalebox{#2}{$\m@th#1\bullet$}}}}}
\makeatother

\usetikzlibrary{calc,decorations.pathmorphing,shapes,cd}

\newcounter{sarrow} % ensures that each node has a unique name

\stackMath
\newcommand\reallywidehat[1]{%
\savestack{\tmpbox}{\stretchto{%
  \scaleto{%
    \scalerel*[\widthof{\ensuremath{#1}}]{\kern.1pt\mathchar"0362\kern.1pt}%
    {\rule{0ex}{\textheight}}%WIDTH-LIMITED CIRCUMFLEX
  }{\textheight}% 
}{2.4ex}}%
\stackon[-6.9pt]{#1}{\tmpbox}%
}
\title[GV invariants of multi-Banana configurations]{Genus Zero Gopakumar-Vafa invariants of multi-Banana configurations}

\begin{document}

\date{\today}

\author{Nina Morishige}
\address{Nina Morishige, Department of Mathematics, The University of British Columbia, Vancouver, BC, V6T 1Z2 Canada}%
\email{nina@math.ubc.ca}%
%\begin{comment}
\begin{abstract}
The multi-Banana configuration $\Lm$ is a local Calabi-Yau threefold of Schoen type. Namely, $\Lm$ is a conifold resolution of $\widehat{I}_v \times_{\DD} \widehat{I}_w$, where  $\widehat{I}_v \to {\DD}$ is an elliptic surface over a formal disc ${\DD}$ with an $I_v$ singulararity on the central fiber. We generalize the technique developed in our earlier paper to compute genus 0 Gopakumar-Vafa invariants of certain fiber curve classes. We illustrate the computation explicitly for $v=1$ and $v=w=2$. The resulting partition function can be expressed in terms of elliptic genera of $\CC^2$, or classical theta functions, respectively.
\end{abstract}
%\end{comment}
\maketitle

\section{Introduction}\label{sec:intro}

\subsection{Background}
Let $X$ be a quasi-projective Calabi-Yau threefold over $\CC$, so that $X$ is smooth and $K_X\cong \OO_X$. Fix a curve class $\beta\in H_2(X)$. Let $M=M^X_{\beta}$ be the moduli space of Simpson semistable \cite{simpson94}, pure, 1-dimensional sheaves $\FF$ with proper support on $X$ with $\text{ch}_2(\FF) = \beta^{\vee}$ and $\chi(\FF) = 1$. The genus 0 Gopakumar-Vafa invariants $n^0_{\beta}(X)$ are defined mathematically by Katz \cite{katz08}:
\begin{definition} 
The genus 0 Gopakumar-Vafa (GV) invariants $n^0_{\beta}(X)$ of $X$ in curve class $\beta$ are defined as the Behrend function weighted Euler characteristics of the moduli space $M^X_{\beta}$.
\begin{equation}
n^0_{\beta}(X) = e(M^X_{\beta},\nu) \coloneqq \sum_{k\in\ZZ}{k\cdot e_{\textsl{top}}(\nu^{-1}(k))}
\label{def:GV}
\end{equation}
where $e_{\textsl{top}}$ is topological Euler characteristic and $\nu:M^X_{\beta}\rightarrow \ZZ$ is Behrend's constructible function \cite{behrend09}.
\end{definition}

In our previous paper \cite{nina19}, we computed the genus 0 Gopakumar-Vafa invariants of the Banana manifold, $\Xb$, a special kind of Schoen threefold, defined as the conifold resolution given by blowing up along the diagonal of the fiber product of a generic rational elliptic surface $S\to\PP^1$ with itself :
\[
\Xb\coloneqq \Bl_{\Delta}(S\times_{\PP^1} S).
\] 
These results were consistent with the computation of the Donaldson-Thomas invariants of $\Xb$ obtained via topological vertex methods by Bryan \cite{bryan19}.

In this paper, we use similar methods as before to obtain the genus 0 Gopakumar-Vafa invariants of certain fiber classes of related local Calabi-Yau threefolds, which we call multi-Banana configurations, and denote by $\Lm$. Our motivation is to study the fiberwise contribution of these configurations, which exist as formal subschemes in special Schoen manifolds, (Section~\ref{subsec:global}). Unlike in our previous paper, even the genus 0 Gopakumar-Vafa invariants associated to these configurations cannot be obtained by other methods at present. Additionaly, the example configurations we study yield partition functions with modular properties that can be expressed succinctly. Our results appear to be compatible with results that appear in the physics literature \cite[Section~3.3]{iqbal13}.

\subsection{The multi-Banana configuration $\Lm$}

The twelve singular fibers $\Fb$ of the regular Banana manifold $\Xb$ are normalizations of the product of $I_1$ singular fibers with themselves, 
\[
\Fb\cong \Bl_{\Delta}(I_1\times I_1)\subset \Xb.
\]
Let $\Fbhat$ be the formal completion of $\Xb$ along $\Fb$. Each $\Fb$ is isomorphic to a non-normal toric variety whose normalization is isomorphic to $\PP^1 \times \PP^1$ blown up at two points on the diagonal. We have $\pi_1(\Fbhat) = \ZZ\times\ZZ$. See \cite[Section~3.1]{nina19} for details.

We define the local Calabi-Yau threefold $\Lm$ as follows:

\begin{definition}\label{def:Lm} 
The multi-Banana $\Fvw{v}{w}$ and the local multi-Banana configuration $\Lvw{v}{w}$ are the \'{e}tale covers of $\Fb$ and $\Fbhat$, respectively,
\[\Fvw{v}{w} \to \Fb,\] 
\[\Lvw{v}{w} \to \Fbhat,\] 
associated to the subgroup $  v\ZZ\times w\ZZ \subset \ZZ\times\ZZ$. 
\end{definition}
We sometimes suppress the decoration and write $\Fm$ and $\Lm$ instead. Observe that $\Fb=\Fvw{1}{1}$ and $\Fbhat=\Lvw{1}{1}$.

The geometry of multi-Banana configurations was studied by Kanazawa and Lau \cite{kanazawa-lau19}. In particular, $\Lvw{v}{w}$ has $vw+2$ curve classes, generated by three families of curves, $\{A_i\}$, $\{B_j\}$, and $\{C_k\}$, see section~\ref{subsec:local}:
\[
\beta \in \sum\limits_{i=0}^{w-1}{ \ZZ [A_i] } \oplus  \sum\limits_{j=0}^{v-1} {\ZZ [B_j]}  \ \oplus \sum\limits_{k=0}^{(v-1)(w-1)} {\ZZ[C_k]}, \qquad  \beta \in H_2(\Lvw{v}{w}).
\]

\subsection{Main results}
In some cases of small $v$ and $w$, the GV invariants have nice formulas. We can express the partition function in terms of $\JphiQ{p}$, the unique weak Jacobi form of weight -2 and index 1,
\begin{align*}
&\JphiQ{p} = p^{-1}(1-p)^2\prod_{m=1}^{\infty}{\frac{(1-Q^{m}p^{-1})^2 (1-Q^{m}p)^2}{(1-Q^m)^4}},\\
Q&=\exp(2\pi i\tau),  \qquad p = \exp(2\pi iz), \qquad  (\tau,z)\in \HH\times\CC.
\end{align*}

and $\mathrm{Ell}_{Q,p}(\CC^2,t)$, the equivariant elliptic genus of $\CC^2$:

\begin{align*}
&\mathrm{Ell}_{Q,p}(\CC^2,t) =\frac{\sqrt{\JphiQ{pt}\JphiQ{p^{-1}t}}}{\JphiQ{t}}.\\ 
%\label{eq:EllipticGenusC2}  
%q&=\exp(2\pi i\tau),  \qquad p = \exp(2\pi iz), \qquad  (\tau,z)\in \HH\times\CC.
\end{align*}

When $v=1$, we have the following.
\begin{theorem} (See Theorem~\ref{thm:1W} for details and notation.)
Fix a curve class $\beta_{(\textbf{a},\,c)}$ in the local multi-Banana $\Lm=\Lvw{1}{w}$: 
\begin{gather*}
\beta_{(\textbf{a},\,c)} = \sum\limits_{i=0}^{w-1}{{a_i}[A_i]} +{c}[C] + [B],\\
\textbf{ a} = (a_0,\ldots,a_{w-1}) \in \ZZ_{\geq 0}^w,\, c  \in \ZZ_{\geq 0}.
\end{gather*}

Then the genus 0 Gopakumar-Vafa invariants $n^0_{\beta_{(\textbf{a},c)}}(\Lvw{1}{w})$ can be expressed as:
\[
\sum_{\textbf{a}, c} {n^0_{\beta_{(\textbf{a},\,{c})}}(\Lm) \textbf{r}^{\textbf{a}} s^{c}} =
s\cdot\JphiQ{s} \sum_{i=0}^{w-1} \quad {\prod\limits_{k=i}^{i+w-2}  \mathrm{Ell}_{Q,s}(\CC^2,R_{i;k}) },
\]
where
\begin{gather*}
%\textbf{r}^{\textbf{a}}  \coloneqq r_0^{a_0}r_1^{a_1}\ldots r_{w-1}^{a_{w-1}},\\
Q\coloneqq \prod\limits_{i=0}^{w-1}{(r_is)},\\
R_{a;b}\coloneqq r_a\cdot r_{a+1}\cdot r_{a+2} \cdots r_b \cdot s^{b-a+1}, \quad a\leq b,\\
r_{k+w} \coloneqq r_k.
\end{gather*} 
\end{theorem}

%
%\subsubsection{Case $2,2$}
In the case of $v=w=2$, the curve classes are naturally labelled as $A_0, A_1, B_0, B_1, C_0, C_1$. We have the following result:

\begin{theorem} (See Theorem~\ref{thm:22} for details and notation.)
Let $v=w=2$, and fix a curve class $\beta_{(\textbf{a},\,\textbf{c})}$ in the local multi-Banana $\Lm=\Lvw{2}{2}$: 
\begin{gather*}
\beta_{(\textbf{a},\,\textbf{c})} = {a_0}[A_0] + {a_1}[A_1] +{c_0}[C_0] + {c_1}[C_1] + [B_0],\\
\textbf{ a} = (a_0,a_1), \textbf{ c} =(c_0,c_1) \in \ZZ_{\geq 0}^2.
\end{gather*}
Then the genus 0 Gopakumar-Vafa invariants $n^0_{\beta_{(\textbf{a},\textbf{c})}}(\Lm)$ are given by the following:
\[
\sum_{a_0, a_1, c_0, c_1} {n^0_{\beta_{(\textbf{a},\textbf{c})}}(\Lm) r_0^{a_0}r_1^{a_1} s_0^{c_0} s_1^{c_1}} =
2\left\{\frac{\JphiQ{r_0}\JphiQ{s_0}\JphiQ{r_1}\JphiQ{s_1} } { \JphiQ{r_0s_0}\JphiQ{r_1s_1} } \right\}^{1/2},
\]
%where $\JphiQ{p}\coloneqq \Jphi{Q}{p}$ is the unique weak Jacobi form of weight -2 and index 1.
where
\begin{align*}
Q&\coloneqq r_0 r_1 s_0 s_1\\
\JphiQ{p}&\coloneqq \Jphi{Q}{p}.
\end{align*}
\end{theorem}

\begin{remark} We note that the appearance of the elliptic genera in the partition function of the multi-Banana suggests a correspondence via geometric engineering \cite{katzklemmvafa97} to partition functions of Yang-Mills gauge theories on surfaces. This viewpoint is discussed further in the previously cited physics literature \cite{iqbal13}.
\end{remark}

%%%%%%%%%%%%%%%%%%%%%%%%%%%%%%%%%%%%%%%%
\subsection{Outline of method} \label{subsec:method}
We recall the method we used in \cite{nina19} to compute the genus 0 GV invariants of $\Xb$. The argument carries over largely unchanged for the local multi-Banana configurations $X=\Lm$, apart from the final combinatorics computation, so we refer the reader to our previous paper for the details of the proofs of the statements in this summary of our method.

We have a $T\coloneqq\tor$ torus action on $\Fm$, given by translation on the smooth locus, and which extends to an action on all of $\Fmhat$. This gives us an action on its coherent sheaves $\Coh(\Fmhat)$ and thus on the moduli space $\MF$. This action preserves the canonical class and is compatible with the symmetric obstruction theory. We can use the motivic nature of the Behrend function weighted Euler characteristic to stratify the moduli space under this group action \cite{behrend09, behrend-fantechi08}. The nontrivial torus orbits make no contribution to $e(\MF,\nu)$, and we can reduce to considering only the $T$-fixed points of the moduli space $(\MF)^T$.

We first count the fixed points of the moduli space. This gives us the naive Euler characteristic, $\naivel$, which we define as the Euler characteristic of the moduli space without the Behrend function weighting:
\[
\naivel \coloneqq e(M^{\Lm}_{\beta}).
\]
Using stability arguments we show that the sheaves in our moduli space have scheme-theoretic support on the multi-Banana surface $\Fm$ \cite[Proposition 12]{nina19}. Thus, for computing $\naivel$, it suffices to count $T$-invariant sheaves of $\Fm$. 

We would like to work on the universal cover of a multi-Banana, $\UFm$, to make the computations easier. This is an infinite type toric surface, whose irreducible components are isomorphic to the blow-up of $\PP^1\times\PP^1$ at two torus fixed points. We give further details of the local geometry in Section~\ref{subsec:local}. The universal cover $\UFm$ is the same as that of the regular Banana fiber $\UFb$ considered in \cite{nina19}. 

In order to relate the sheaves of $\UFm$ with those of $\Fm$, we introduce another torus action, which we denote by $P\coloneqq\tor$. This $P$ action on $\Coh(\Fm)$ is defined by tensoring with degree 0 line bundles of $\Fm$, as $\Pic^0(\Fm)\cong\tor$ \cite[Section 4]{nina19}.  Again, the Euler characteristic contribution can be computed on orbits of the action, and it then suffices to consider only sheaves invariant under the two $\tor$ actions, $T$ and $P$. The $T$ torus action also lifts to give an action on the universal cover and its sheaves. 

Any sheaf fixed under the $P$ action pulls back to an equivariant sheaf on $\UFm$ \cite[Proposition~22]{nina19}. This equivariant sheaf contains a distinguished subsheaf which pushes forward to the original sheaf, and is unique up to deck transformations. Moreover stability, Euler characteristic 1, and invariance under the $T$ torus action is preserved in this correspondence. 

The requirements of stability and Euler characteristic equal to 1 then puts restrictions on the allowed invariant stable sheaves. If we further specify that the curve class $\beta$ has degree exactly 1 in one of the curve families of $\Fvw{v}{w}$, then all the $T$ and $P$ fixed sheaves in our moduli space correspond to structure sheaves of possibly non-reduced curves on $\UFm$ \cite[Proposition~23]{nina19}. The multiplicity of each component is constrained \cite[Proposition~31]{nina19} by a condition, which is equivalent to requiring that the partition given by multiplicities of successive rational components from the fixed central degree 1 curve has a conjugate partition with odd parts that are distinct. We give the specific details of this condition in Section~\ref{sec:22}.

This count of the number of fixed points of the moduli space gives the naive Euler characteristic, $\naivel$. However, the Behrend function weighting amounts to a sign $(-1)^{\deg\beta}$, which depends on the total degree of the curve \cite[Proposition~39]{nina19}, and this can be incorporated into the partition function.

We note that our technique is limited to computing invariants associated to fiber class curves such that the degree of one of these families is fixed to be 1. We do not yet know how to extend the technique to arbitrary degrees.

Our method allows us to calculate the partition function for $\Lvw{v}{w}$ in the general case, for arbitrary $v$ and $w$. However, as $v$ and $w$ increase, there will be unavoidable linear relations among the curve classes, even after fixing the degree of one curve type to be 1. We will only present in detail the $2\times 2$ case (Section~\ref{sec:22}) and the $1\times w$ case (Section~\ref{sec:1w}) as they illustrate the ideas sufficiently without the notation becoming burdensome.

%%%%%%%%%%%%%%%%%%%%%%%%%%%%%%%%%%%%%%%%%%%
%\input{./section1}
%\newpage
\section{Geometry}
\label{sec:geometry}
In this section, we give two examples of multi-Banana configurations $\Lm$ that exist as formal neighborhoods of surfaces inside compact Calabi-Yau threefolds. We then discuss some of the local geometry of multi-Banana configurations needed for the following sections.

\subsection{Global geometry} 
\label{subsec:global}
\begin{definition}
A multi-Banana manifold $\Xm$ is a smooth Calabi-Yau threefold which is a conifold resolution of the fiber product of two rational elliptic surfaces, and such that the formal neighborhood of each singular fiber is a multi-Banana configuration.
\end{definition}

\begin{example}
\label{ex:2x2}
Let $S\stackrel{\pi}{\to} \PP^1$ be a rational elliptic surface with singular fibers consisting of four $I_1$ and four $I_2$ singular fibers. Suppose $S$ has a $2$-torsion section. This induces an order $2$ automorphism $\phi_2$ that interchanges the nodes of each of the $I_2$ fibers. 

We can then form the fiber product of $S$ with itself, $S\times_{\PP^1} S$. In order to get a conifold resolution, we blow up the generalized diagonal $\widetilde{\Delta}$, consisting of the diagonal $\Delta$, as well as all its translates by iterations of $\phi_2$,
\begin{equation*}
\widetilde{\Delta}\coloneqq (\phi_2^{(i)} \times \phi_2^{(j)}) \Delta,\quad 0\leq i, j <2.
\end{equation*}

We will call this multi-Banana manifold $\Xd$.
\begin{equation*}
\Xd \coloneqq \Bl_{\widetilde{\Delta}}(S\times_{\PP^1} S)
\end{equation*}
In this case, the multi-Banana contains four $\Fvw{2}{2}$ configurations, and four ordinary Banana fibers $\Fb$.

\end{example}

Instead of taking the fiber product of $S$ with itself, we can also do the following construction. 
\begin{example}
\label{ex:1x5}
Let $S\stackrel{\pi}{\to} \PP^1$ be a rational elliptic surface with two $I_1$ and two $I_5$ singular fibers. Then $S$ has a 5-torsion section, which induces an order $5$ automorphism $\phi_5$ which acts on each $I_5$ fiber by cycling the nodes.

Now, let us take the quotient of $S$ by the action of $\phi_5$, and let $S'$ be the resolution of the quotient: 
$$S'\coloneqq \text{Res\,}(S/\phi_5).$$
Notice that by construction, $S'\stackrel{\pi'}{\to}\PP^1$ is another rational elliptic surface with singular fibers over the same base points:
\begin{align*}
\PP^1_{\mathsc{sing}}\coloneqq&\{p \in \PP^1\vert \pi^{-1}(p) \subset S \text{ singular} \} \\
=& \{p \in \PP^1\vert {\pi'}^{-1}(p) \subset S' \text{ singular} \}. 
\end{align*}
We also have that the smooth fibers of $S$ and $S'$ are isogenous: 
\begin{gather*}
\phi_5|_p:\pi^{-1}(p) \to {\pi'}\,^{-1}(p), \\
\qquad \phi_5|_p \text{ is an isogeny, } \forall p \in \PP^1 \backslash\PP^1_{\mathsc{sing}}.
\end{gather*}
In this case, there is a conifold resolution of the fiber product, $S\times_{\PP^1} S'$. From the construction, we have a rational map of schemes over $\PP^1$, 
\[
S \dashrightarrow S'=\text{Res\,}(S/\phi_k),
\]
so we get a graph 
\[
\overline{\Gamma}\subset S\times_{\PP^1} S'.
\]
Then the conifold resolution is given by blowing up this graph $\overline{\Gamma}$. We will call this multi-Banana threefold $\Xq$:
\[
\Xq \coloneqq \Bl_{\overline{\Gamma}} (S\times_{\PP^1} S').
\]
The multi-Banana manifold $\Xq$ is a rigid Calabi-Yau threefold and contains two $\Fvw{1}{5}$ multi-Banana configurations, and two $\Fvw{5}{1}$ multi-Banana configurations.
\end{example}
%
%
%%%%%%%%%%%%%%%%%%%%%%%%%%%%%%%%%%%%%%%%%%%%
\subsection{Local geometry}\label{subsec:local}
We now examine the local geometry of the multi-Banana configurations in more detail and establish some notation we will need later.

We recall the construction from \cite{kanazawa-lau19} and relate it to our discussion.

Let $L$ be a tiling of the plane $(x,y,1)\subset \RR^3$ given by: 
%\begin{equation*}
\begin{gather*}
%\begin{aligned}
\{x=m, z=1\} \cup \{y=n, z=1\} \cup \{y-x=r,  z=1\}, \\
 m,n,r\in \ZZ,\\  (x,y,z)\in \RR^3. 
%\end{aligned}
\end{gather*}
%\end{equation*}
  
\begin{figure}[th!]
  \centering
	\scalebox{.5}{

\begin{tikzpicture}[y=0.80pt, x=0.80pt, yscale=-1.000000, xscale=1.000000, inner sep=0pt, outer sep=0pt]
\begin{scope}[shift={(0,0)}]
  \path[fill=black,line join=miter,line cap=butt,line width=0.800pt]
    (0.0000,0.0000) node[above right] (flowRoot4196) {};
      \path[draw=black,line join=miter,line cap=butt,miter limit=4.00,even odd
        rule,line width=1.200pt] (102.8641,70.6532) -- (410.8945,71.4646) --
        (410.8945,71.4646);
      \path[draw=black,line join=miter,line cap=butt,miter limit=4.00,even odd
        rule,line width=1.200pt] (169.5110,52.5395) -- (480.3115,52.9452) --
        (480.3115,52.9452);
      \path[draw=black,line join=miter,line cap=butt,miter limit=4.00,even odd
        rule,line width=1.200pt] (26.5749,91.3921) -- (337.2085,91.7977) --
        (337.2085,91.7977);
      \path[draw=black,line join=miter,line cap=butt,miter limit=4.00,even odd
        rule,line width=1.200pt] (114.8988,93.3662) -- (279.7972,48.3332);
      \path[draw=black,line join=miter,line cap=butt,miter limit=4.00,even odd
        rule,line width=1.200pt] (208.5741,93.3648) -- (373.2696,49.3460);
      \path[draw=black,line join=miter,line cap=butt,miter limit=4.00,even odd
        rule,line width=1.200pt] (307.7416,93.7659) -- (475.0740,48.9358);
      \path[draw=black,line join=miter,line cap=butt,miter limit=4.00,even odd
        rule,line width=1.200pt] (106.2611,66.7448) -- (127.4439,98.8574);
      \path[draw=black,line join=miter,line cap=butt,miter limit=4.00,even odd
        rule,line width=1.200pt] (176.8772,50.8156) -- (218.3697,94.6056);
      \path[draw=black,line join=miter,line cap=butt,miter limit=4.00,even odd
        rule,line width=1.200pt] (251.2464,41.2304) -- (328.1361,102.3769);
      \path[draw=black,line join=miter,line cap=butt,miter limit=4.00,even odd
        rule,line width=1.200pt] (348.5892,45.1502) -- (398.8461,78.5130);
      \path[draw=black,line join=miter,line cap=butt,miter limit=4.00,even odd
        rule,line width=1.200pt] (447.2707,41.4132) -- (472.4974,61.5030);
      \path[draw=black,line join=miter,line cap=butt,miter limit=4.00,even odd
        rule,line width=0.800pt] (196.5251,16.8382) -- (195.7138,259.8516);
      \path[draw=black,line join=miter,line cap=butt,miter limit=4.00,even odd
        rule,line width=0.800pt] (47.7591,181.5541) -- (394.2752,181.5541);
      \path[draw=black,line join=miter,line cap=butt,miter limit=4.00,even odd
        rule,line width=0.800pt] (59.5504,217.7284) -- (319.1645,148.2695);
    \path[draw=black,line join=miter,line cap=butt,miter limit=4.00,even odd
      rule,line width=0.600pt] (195.7014,181.3331) -- (122.2716,91.7768);
    \path[draw=black,line join=miter,line cap=butt,miter limit=4.00,even odd
      rule,line width=0.600pt] (195.6490,180.8086) -- (215.5274,91.6334);
    \path[draw=black,line join=miter,line cap=butt,miter limit=4.00,even odd
      rule,line width=0.600pt] (196.0777,180.3183) -- (314.6478,91.8453);
    \path[draw=black,line join=miter,line cap=butt,miter limit=4.00,even odd
      rule,line width=0.600pt] (195.6597,181.3573) -- (359.9582,52.8154);
    \path[draw=black,line join=miter,line cap=butt,miter limit=4.00,even odd
      rule,line width=0.600pt] (195.8371,181.0154) -- (289.1078,71.7498);
    \path[draw=black,line join=miter,line cap=butt,miter limit=4.00,even odd
      rule,line width=0.600pt] (388.2839,71.4454) -- (195.7599,181.0058);
    \path[draw=black,line join=miter,line cap=butt,miter limit=4.00,even odd
      rule,line width=0.600pt] (266.0245,53.5603) -- (195.7910,181.1275);
    \path[draw=black,line join=miter,line cap=butt,miter limit=4.00,even odd
      rule,line width=0.600pt] (195.1451,181.2890) -- (34.7967,91.1590);
    \path[draw=black,line join=miter,line cap=butt,miter limit=4.00,even odd
      rule,line width=0.600pt] (461.4305,52.5124) -- (195.6085,181.2491);
    \path[draw=black,line join=miter,line cap=butt,miter limit=4.00,even odd
      rule,line width=1.200pt] (26.9818,93.6945) -- (191.8802,48.6616);
    \path[draw=black,line join=miter,line cap=butt,miter limit=4.00,even odd
      rule,line width=1.200pt] (32.2834,87.1724) -- (38.6942,96.2594);
    \path[draw=black,line join=miter,line cap=butt,miter limit=4.00,even odd
      rule,line width=0.600pt] (195.5719,181.1696) -- (108.6597,71.2463);
    \path[draw=black,line join=miter,line cap=butt,miter limit=4.00,even odd
      rule,line width=0.600pt] (195.6181,180.8026) -- (177.8392,52.8070);
    \begin{scope}[cm={{-0.92337,0.38391,-0.38391,-0.92337,(519.98043,328.00729)}}]
      \path[draw=black,fill=black,miter limit=1.00,even odd rule,line width=0.279pt]
        (333.1611,333.7799)arc(0.000:89.773:1.743)arc(89.773:179.547:1.743)arc(179.547:269.321:1.743)arc(269.320:359.094:1.743);
      \path[draw=black,fill=black,miter limit=1.00,even odd rule,line width=0.279pt]
        (323.1611,333.7799)arc(0.000:89.773:1.743)arc(89.773:179.547:1.743)arc(179.547:269.321:1.743)arc(269.320:359.094:1.743);
      \path[draw=black,fill=black,miter limit=1.00,even odd rule,line width=0.279pt]
        (343.1611,333.7799)arc(0.000:89.773:1.743)arc(89.773:179.547:1.743)arc(179.547:269.321:1.743)arc(269.320:359.094:1.743);
    \end{scope}
    \path[draw=black,fill=black,miter limit=1.00,even odd rule,line width=0.279pt]
      (470.2927,73.9342)arc(-0.000:89.774:1.743)arc(89.773:179.547:1.743)arc(179.547:269.321:1.743)arc(269.320:359.094:1.743);
    \path[draw=black,fill=black,miter limit=1.00,even odd rule,line width=0.279pt]
      (460.2928,73.9342)arc(-0.000:89.774:1.743)arc(89.773:179.547:1.743)arc(179.547:269.321:1.743)arc(269.320:359.094:1.743);
    \path[draw=black,fill=black,miter limit=1.00,even odd rule,line width=0.279pt]
      (480.2927,73.9342)arc(-0.000:89.774:1.743)arc(89.773:179.547:1.743)arc(179.547:269.321:1.743)arc(269.320:359.094:1.743);
    \path[draw=black,fill=black,miter limit=1.00,even odd rule,line width=0.279pt]
      (30.5547,64.2591)arc(-0.000:89.774:1.743)arc(89.774:179.547:1.743)arc(179.547:269.320:1.743)arc(269.320:359.094:1.743);
    \path[draw=black,fill=black,miter limit=1.00,even odd rule,line width=0.279pt]
      (20.5548,64.2591)arc(-0.000:89.773:1.743)arc(89.773:179.547:1.743)arc(179.547:269.321:1.743)arc(269.321:359.094:1.743);
    \path[draw=black,fill=black,miter limit=1.00,even odd rule,line width=0.279pt]
      (40.5547,64.2591)arc(-0.000:89.774:1.743)arc(89.774:179.547:1.743)arc(179.547:269.320:1.743)arc(269.320:359.094:1.743);
    \begin{scope}[cm={{-0.92337,0.38391,-0.38391,-0.92337,(772.67419,213.06684)}}]
      \path[draw=black,fill=black,miter limit=1.00,even odd rule,line width=0.279pt]
        (333.1611,333.7799)arc(0.000:89.773:1.743)arc(89.773:179.547:1.743)arc(179.547:269.321:1.743)arc(269.320:359.094:1.743);
      \path[draw=black,fill=black,miter limit=1.00,even odd rule,line width=0.279pt]
        (323.1611,333.7799)arc(0.000:89.773:1.743)arc(89.773:179.547:1.743)arc(179.547:269.321:1.743)arc(269.320:359.094:1.743);
      \path[draw=black,fill=black,miter limit=1.00,even odd rule,line width=0.279pt]
        (343.1611,333.7799)arc(0.000:89.773:1.743)arc(89.773:179.547:1.743)arc(179.547:269.321:1.743)arc(269.320:359.094:1.743);
    \end{scope}
\end{scope}

\end{tikzpicture}
	
	}
  \caption{The fan of $\mathcal{A}$.}
	\label{fig:infinitetilecone}
\end{figure}
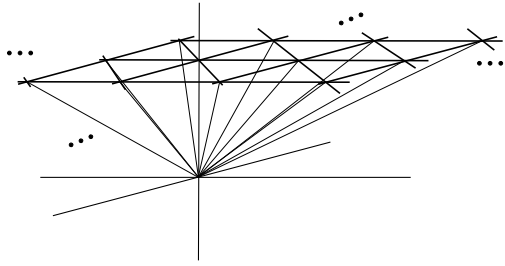

Let $\mathcal{A}$ be the non-finite type toric threefold whose fan consists of all the cones over the proper faces of $L$ (Figure~\ref{fig:infinitetilecone}). Let $\UFm$ be the universal cover of $\Fm$, and $\UFmhat$ the universal cover of $\Fmhat$, 
\begin{align*}
\UFm&\xrightarrow{~pr~}\Fm  \\
\UFmhat&\xrightarrow{~pr~}\Fmhat.
\end{align*}
Then $\UFm\subset \mathcal{A}$ is the union of the toric divisors of $\mathcal{A}$, and $\UFmhat$ is the formal completion of $\mathcal{A}$ along $\UFm$.

We have an action of $G=v\ZZ \times w\ZZ$, $v, w \in \ZZ_{\geq 0}$, on $L\subset\RR^3$ by translation:
\[
(v,w)\cdot (x,y,1) = (x+v, y+w,1).
\]
which induces an automorphism $\psi_G:\mathcal{A}\to \mathcal{A}$ and also on $\UFmhat$. These are then the deck transformations of the universal cover of the local multi-Banana configuration $\Lm$:
\begin{center}
\[
\UFmhat \to \UFmhat/\psi_G \cong \Lm.
\]
\end{center}

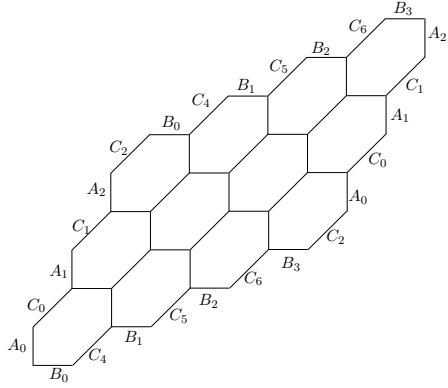
\begin{figure}[th!]
  \centering
	\scalebox{.7}{

\begin{tikzpicture}[y=0.80pt, x=0.80pt, yscale=-1.000000, xscale=1.000000, inner sep=0pt, outer sep=0pt]

\node at (47,275) {\scalebox{.8}{$A_0$}};
\node at (75,225) {\scalebox{.8}{$A_1$}};
\node at (100,170) {\scalebox{.8}{$A_2$}};
\node at (275,175) {\scalebox{.8}{$A_0$}};
\node at (303,120) {\scalebox{.8}{$A_1$}};
\node at (328,65) {\scalebox{.8}{$A_2$}};

\node at (75,295) {\scalebox{.8}{$B_0$}};
\node at (125,270) {\scalebox{.8}{$B_1$}};
\node at (175,245) {\scalebox{.8}{$B_2$}};
\node at (230,220) {\scalebox{.8}{$B_3$}};
\node at (150,125) {\scalebox{.8}{$B_0$}};
\node at (200,100) {\scalebox{.8}{$B_1$}};
\node at (250,75) {\scalebox{.8}{$B_2$}};
\node at (305,48) {\scalebox{.8}{$B_3$}};

\node at (60,248) {\scalebox{.8}{$C_0$}};
\node at (90,195) {\scalebox{.8}{$C_1$}};
\node at (115,142) {\scalebox{.8}{$C_2$}};
\node at (288,150) {\scalebox{.8}{$C_0$}};
\node at (313,100) {\scalebox{.8}{$C_1$}};
\node at (260,202) {\scalebox{.8}{$C_2$}};

\node at (101,282) {\scalebox{.8}{$C_4$}};
\node at (155,255) {\scalebox{.8}{$C_5$}};
\node at (205,230) {\scalebox{.8}{$C_6$}};
\node at (170,110) {\scalebox{.8}{$C_4$}};
\node at (222,87) {\scalebox{.8}{$C_5$}};
\node at (275,60) {\scalebox{.8}{$C_6$}};

  \begin{scope}[cm={{0.3239,0.0,0.0,0.3239,(-101.67469,46.60915)}}]
  \end{scope}
  \begin{scope}[cm={{0.3239,0.0,0.0,0.3239,(-75.69872,-5.19936)}}]
  \end{scope}
  \begin{scope}[cm={{0.3239,0.0,0.0,0.3239,(-49.72276,-57.00787)}}]
  \end{scope}
  \begin{scope}[cm={{0.59536,0.0,0.0,0.59536,(-52.69554,-126.12323)}}]
  \end{scope}
  \begin{scope}[cm={{0.3239,0.0,0.0,0.3239,(-101.14992,124.19174)}}]
  \end{scope}
    \path[draw=black,line join=round,line cap=round,even odd rule,line
      width=0.259pt] (83.5379,288.4620) -- (57.6257,288.4620) -- (57.6257,262.6309)
      -- (83.5379,236.7187);
    \path[draw=black,line join=round,line cap=round,even odd rule,line
      width=0.259pt] (109.5138,236.6535) -- (83.6017,236.6535) -- (83.6017,210.8223)
      -- (109.5138,184.9102);
    \path[draw=black,line join=round,line cap=round,even odd rule,line
      width=0.259pt] (135.4898,184.8450) -- (109.5777,184.8450) --
      (109.5777,159.0138) -- (135.4898,133.1017);
    \path[draw=black,line join=round,line cap=round,even odd rule,line
      width=0.259pt] (161.4658,133.0365) -- (135.5536,133.0365);
    \path[draw=black,line join=round,line cap=round,even odd rule,line
      width=0.259pt] (84.1265,288.4049) -- (110.0386,262.4928);
    \path[draw=black,line join=round,line cap=round,even odd rule,line
      width=0.259pt] (136.0146,262.4276) -- (110.1024,262.4276) --
      (110.1024,236.5964) -- (136.0146,210.6843);
    \path[draw=black,line join=round,line cap=round,even odd rule,line
      width=0.259pt] (161.9905,210.6191) -- (136.0784,210.6191) --
      (136.0784,184.7879) -- (161.9905,158.8758);
    \path[draw=black,line join=round,line cap=round,even odd rule,line
      width=0.259pt] (187.9665,158.8106) -- (162.0544,158.8106) --
      (162.0544,132.9794) -- (187.9665,107.0673);
    \path[draw=black,line join=round,line cap=round,even odd rule,line
      width=0.259pt] (213.9425,107.0021) -- (188.0303,107.0021);
    \path[draw=black,line join=round,line cap=round,even odd rule,line
      width=0.259pt] (136.6032,262.3705) -- (162.5153,236.4584);
    \path[draw=black,line join=round,line cap=round,even odd rule,line
      width=0.259pt] (188.4913,236.3932) -- (162.5791,236.3932) --
      (162.5791,210.5620) -- (188.4913,184.6499);
    \path[draw=black,line join=round,line cap=round,even odd rule,line
      width=0.259pt] (214.4672,184.5847) -- (188.5551,184.5847) --
      (188.5551,158.7535) -- (214.4672,132.8414);
    \path[draw=black,line join=round,line cap=round,even odd rule,line
      width=0.259pt] (240.4432,132.7761) -- (214.5310,132.7761) --
      (214.5310,106.9450) -- (240.4432,81.0328);
    \path[draw=black,line join=round,line cap=round,even odd rule,line
      width=0.259pt] (266.4192,80.9676) -- (240.5070,80.9676);
    \path[draw=black,line join=round,line cap=round,even odd rule,line
      width=0.259pt] (189.0798,236.3361) -- (214.9920,210.4239);
    \path[draw=black,line join=round,line cap=round,even odd rule,line
      width=0.259pt] (240.9680,210.3587) -- (215.0558,210.3587) --
      (215.0558,184.5276) -- (240.9680,158.6154);
    \path[draw=black,line join=round,line cap=round,even odd rule,line
      width=0.259pt] (266.9439,158.5502) -- (241.0318,158.5502) --
      (241.0318,132.7191) -- (266.9439,106.8069);
    \path[draw=black,line join=round,line cap=round,even odd rule,line
      width=0.259pt] (292.9199,106.7417) -- (267.0077,106.7417) --
      (267.0077,80.9106) -- (292.9199,54.9984);
    \path[draw=black,line join=round,line cap=round,even odd rule,line
      width=0.259pt] (318.8958,54.9332) -- (292.9837,54.9332);
    \path[draw=black,line join=round,line cap=round,even odd rule,line
      width=0.259pt] (241.5565,210.3017) -- (267.4687,184.3895);
    \path[draw=black,line join=round,line cap=round,even odd rule,line
      width=0.259pt] (267.5325,184.3243) -- (267.5325,158.4931) --
      (293.4446,132.5810);
    \path[draw=black,line join=round,line cap=round,even odd rule,line
      width=0.259pt] (293.5085,132.5158) -- (293.5085,106.6846) --
      (319.4206,80.7725);
    \path[draw=black,line join=round,line cap=round,even odd rule,line
      width=0.259pt] (319.4844,80.7073) -- (319.4844,54.8761);
  \begin{scope}[cm={{0.59536,0.0,0.0,0.59536,(287.31693,62.89966)}}]
  \end{scope}
  \begin{scope}[shift={(428.79521,-270.2712)}]
  \end{scope}

\end{tikzpicture}

	}
  \caption{$\Fm$ in the case $v=3$ and $w=4$. Here, the top boundary curves are identified with those along the bottom, and also the left edge with the right edge.}
	\label{fig:fmban}
\end{figure}

We denote by $\Xi$ the irreducible surface which is the momentum polytope of $\PP^1\times \PP^1$ blown up at 2 points, and drawn as a hexagon in our diagrams: 
\[
\Xi \coloneqq \Bl_{p_1,p_2}(\PP^1\times\PP^1), \quad p_1,p_2\in \PP^1\times\PP^1. 
\]
Then the momentum polytope of $\UFm$ can be represented as a hexagonal tiling of the plane, and the momentum polytope of $\Fm$ is given by $v\times w$ hexagons glued together along their toric boundary, as depicted in the example of Figure~\ref{fig:fmban}. 

The irreducible components of the torus fixed curves in $\Fm$ fall into three families of rational curves. Two of these families, $\{A_i\}$ and $\{B_j\}$, are proper transforms of the rational curves from the $I_v$ and $I_w$ singular fibers in $\Fm$, respectively, and one family, $\{C_k\}$, are the exceptional curves of the conifold resolution.

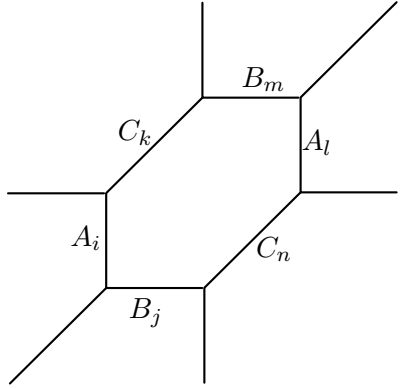
\begin{figure}[th!]
  \centering
	\scalebox{1}{

\begin{tikzpicture}[y=0.80pt, x=0.80pt, yscale=-1.000000, xscale=1.000000, inner sep=0pt, outer sep=0pt]

\node at (57,130) {\scalebox{1}{$A_i$}};
\node at (85,165) {\scalebox{1}{$B_j$}};
\node at (80,80) {\scalebox{1}{$C_k$}};
\node at (165,85) {\scalebox{1}{$A_l$}};
\node at (140,55) {\scalebox{1}{$B_m$}};
\node at (145,135) {\scalebox{1}{$C_n$}};

  \begin{scope}[cm={{0.56129,0.0,0.0,0.56129,(-209.57956,-265.64876)}}]
    \path[draw=black,line join=round,line cap=round,even odd rule,line
      width=0.800pt] (490.0000,667.1122) -- (410.0000,667.1122);
  \end{scope}
  \begin{scope}[cm={{0.56129,0.0,0.0,0.56129,(-160.19283,5.33479)}}]
  \end{scope}
  \begin{scope}[cm={{0.56129,0.0,0.0,0.56129,(-208.67019,-131.20558)}}]
    \path[draw=black,line join=round,line cap=round,even odd rule,line
      width=0.800pt] (410.0000,587.3622) -- (490.0000,507.3622);
  \end{scope}
  \begin{scope}[cm={{0.56129,0.0,0.0,0.56129,(-163.65634,-220.98475)}}]
    \path[draw=black,line join=round,line cap=round,even odd rule,line
      width=0.800pt] (490.0000,667.1122) -- (410.0000,667.1122) --
      (410.0000,587.3622) -- (490.0000,507.3622);
  \end{scope}
  \begin{scope}[cm={{0.56129,0.0,0.0,0.56129,(-118.6425,-310.76392)}}]
    \path[draw=black,line join=round,line cap=round,even odd rule,line
      width=0.800pt] (490.0000,667.1122) -- (410.0000,667.1122) --
      (410.0000,587.3622);
  \end{scope}
  \begin{scope}[cm={{0.56129,0.0,0.0,0.56129,(-117.73313,-176.32075)}}]
    \path[draw=black,line join=round,line cap=round,even odd rule,line
      width=0.800pt] (410.0000,667.1122) -- (410.0000,587.3622) --
      (490.0000,507.3622);
  \end{scope}
  \begin{scope}[cm={{0.56129,0.0,0.0,0.56129,(-72.71927,-266.09991)}}]
    \path[draw=black,line join=round,line cap=round,even odd rule,line
      width=0.800pt] (490.0000,667.1122) -- (410.0000,667.1122) --
      (410.0000,587.3622) -- (490.0000,507.3622);
  \end{scope}
  \begin{scope}[cm={{0.56129,0.0,0.0,0.56129,(-22.42317,139.32681)}}]
  \end{scope}
  \begin{scope}[cm={{0.77267,0.0,0.0,0.77267,(329.27492,-242.0087)}}]
  \end{scope}
  \begin{scope}[cm={{0.77267,0.0,0.0,0.77267,(268.56067,66.65595)}}]
  \end{scope}

\end{tikzpicture}
	}
  \caption{Curve labels for hexagon $\Xi$}
	\label{fig:1x1}
\end{figure}

We will draw these curves oriented as shown in Figure~\ref{fig:1x1}, so the vertical curves are in the $A$ family, the  horizontal curves are $B$ family, and the diagonal curves are from the $C$ family.

When there is no confusion, we will also label irreducible components of the lifts to the universal cover of these torus invariant curves with the curve class of their projection to $\Fm$. That is, an irreducible component of $pr^{-1}(A_i)\subset\UFm$ will also be referred to as $A_i$ in the universal cover.

As each hexagon surface $\Xi$ is the momentum polytope of $\PP^1\times \PP^1$ blown up at 2 points, their 6 boundary divisors have 2 relations. For example, in the labeled Figure~\ref{fig:1x1}, we have:
\begin{equation}
\begin{aligned}
A_i+C_k&=A_l+C_n,\\
B_m+C_k&=B_j+C_n.
\end{aligned}
\label{eq:hex}
\end{equation}
from the equivalent ways to express the total transform of the rulings in each $\PP^1$ factor. 

A priori, there are $3vw$ torus invariant irreducible curves in $\Fm$, but these satisfy standard hexagon relations (Eq.~\ref{eq:hex}), so it is possible to choose a basis of $vw + 2$ curves, consisting of $v \times A$ curves, $w\times B$ curves, and $(v-1)(w-1)+1 \times C$ curves. For the small examples we consider, we will index the curves in each family in a simple way. It is possible to use a systematic choice of generators and labels for the general case \cite[Section~5.1]{kanazawa-lau19}, but it would be notationally cumbersome for these examples, so we do not present that here.

As remarked in the introduction, our technique is limited to considering only curves where we restrict the degree of one family of curves to be exactly 1. For concreteness, we will assume our curves are of class
\begin{equation}
\beta=\sum{a_i[A_i]} +[B_0]+ \sum{c_k[C_k]}.
\label{eq:singleb}
\end{equation}

The hexagonal tiling from $\UFm$ possesses a $v\ZZ\times w\ZZ$ periodicity from the deck transformations. In order to get rid of the ambiguity from the deck transformations, we will assume we have fixed a choice of fundamental domain $D$, and any curve we consider has its unique irreducible component covering $B_0$ inside $\overline{D}$. In other words, we will require that $T$-torus invariant curves $\CCC\subset\UFm$ with $[pr(\CCC)]=\beta$ also satisfy $\CCC\cap pr^{-1}(B_0)\subset \overline{D}$.
 
From the arguments given in Subsection~\ref{subsec:method}, in order to compute the naive Euler characteristic $\naivel$, it suffices to count all configurations of possibly non-reduced $T$-torus invariant curves covering $\beta$ on the universal cover $\UFm$, subject to the constraint that the partition given by multiplicities of successive rational components of each tree emanating from $B_0$ has a conjugate partition which has odd parts that are distinct. In section~\ref{sec:22} and \ref{sec:1w}, we will illustrate this count in two specific cases, namely when the fundamental domain consists of $2\times 2$ hexagons, and also the case of $1\times w$ hexagons. These configurations exist, for example, in $\Xd$, and $\Xq$, respectively, as described in the previous Section~\ref{subsec:global}.

%%%%%%%%%%%%%%%%%%%%%%%%%%%%%%%%%%%%%%%%%%%
%\input{./section2}
\section{Notation and conventions}
We gather in this section the conventions we use for product and sum expansions for Jacobi forms and elliptic genera.  

Recall, the weak Jacobi form $\Jphiqp$ of weight -2 index 1 is defined as 
\[
\Jphiqp = p^{-1}(1-p)^2\prod_{m=1}^{\infty}{\frac{(1-q^{m}p^{-1})^2 (1-q^{m}p)^2}{(1-q^m)^4}},
\]
%\[q=\exp(2\pi i\tau),  \qquad p = \exp(2\pi iz), \qquad  (\tau,z)\in \HH\times\CC.\]

The Jacobi theta function $\theta_1(q,p)$ function is given as
\begin{align*}
\theta_1(q,p) 
&= -\sum{q^{\frac{k^2}{2\ }}(-p)^k}\\
&=-iq^{\frac{1}{8}} p^{-\frac{1}{2}}\prod_{m=1}^{\infty}{(1-q^m)(1-q^{m-1}p)(1-q^{m}p^{-1})}
\end{align*}

and the Dedkind $\eta$ function is
\[
\eta(q)= q^{\frac{1}{24}} \prod_{m=1}^{\infty} (1-q^m).
\]
Here, 
\[q=\exp(2\pi i\tau),  \qquad p = \exp(2\pi iz), \qquad  (\tau,z)\in \HH\times\CC.\]

Since the first variable will be constant within our partition functions, we will use the shortened notation,
\begin{align*}
\JphiQ{p}&\coloneqq \Jphi{Q}{p} \\
\thetaQ{p}&\coloneqq \theta_{1}(Q,p) \\
\etaQ&\coloneqq \eta(Q) 
\end{align*}

Treating these expressions as formal power series, it is easy to verify the identities :
\begin{equation}
\begin{aligned}
\sqrt{\JphiQ{p}} &= \frac{i\thetaQ{p}}{\etaQ^3} ,\\
\sqrt{p\,\JphiQ{p}} &= \sqrt{\frac{Q}{p}\,\JphiQ{\frac{Q}{p}}} ,\\
\sqrt{\JphiQ{p}} &= -\sqrt{\JphiQ{p^{-1}}} .
\end{aligned}
\label{eq:Jphi-identities}
\end{equation}

Suppose $M$ is a non-compact complex manifold of dimension $d$ with a $\CC^*$ action with isolated fixed points $\{x\}$ of tangent weights $k_i$. We define the equivariant elliptic genus of $M$ to be: 
\[
\mathrm{Ell}_{q,y}(M,t)=\sum_{x\in M^{\CC^*}}{\prod_{j=1}^{d}y^{-\frac{1}{2}} \prod_{m=1}^{\infty}{\frac{(1-q^{m-1}yt^{-k_j(x)})(1-q^{m}p^{-1}t^{k_j(x)})}{(1-q^{m-1}t^{-k_j(x)})(1-q^{m}t^{k_j(x)}) }}}.
\]
In particular (\cite[Theorem 12]{waelder08}), we have:
\begin{align}
\mathrm{Ell}_{q,y}(\CC^2,t) &= \frac{\theta_1(q,yt)\theta_1(q,yt^{-1})}{\theta_1(q,t)\theta_1(q,t^{-1})}\nonumber \\
&=\frac{\sqrt{\Jphi{q}{yt}\Jphi{q}{y^{-1}t}}}{\Jphi{q}{t}}\nonumber\\
&= \frac{\sqrt{\JphiQ{yt}\JphiQ{y^{-1}t}}}{\JphiQ{t}} .
\label{eq:EllipticGenusC2}  
\end{align}

%%%%%%%%%%%%%%%%%%%%%%%%%%%%%%%%%%%%%%%%%%%
%\input{./section3}

\section{Case $2\times 2$}\label{sec:22}
In this section, we study the case of $\Ld\coloneqq\Lvw{2}{2}$, and we use this example to illustrate in detail how the method described in the Section~\ref{subsec:method} leads to the computation of the Gopakumar-Vafa invariants. 

\subsection{$T$-Torus fixed curves on $\Ld$} We first fix a choice of a fundamental domain $D$ in $U(\Fd)$ so that there is no ambiguity in our counts due to deck transformations.

We label the curves of the $2\times 2$ hexagons of the momentum polytope of the fundamental domain in $U(\Fd)$ with our convention explained in Section~\ref{subsec:local}. The vertical curves cover curves in the $A$ family, the  horizontal curves those in the $B$ family, and the diagonal curves cover the $C$ family as shown in Figure~\ref{fig:2hexagon}. There is a $2\ZZ\times 2\ZZ$ periodicity of this fundamental domain in the universal cover.

\begin{figure}[th!]
  \centering
	\scalebox{.5}{

\begin{tikzpicture}[y=0.80pt, x=0.80pt, yscale=-1.000000, xscale=1.000000, inner sep=0pt, outer sep=0pt]

\node at (50,385) {\scalebox{1.2}{$A_0$}};
\node at (120,245) {\scalebox{1.2}{$A_1$}};
\node at (190,320) {\scalebox{1.2}{$A_2$}};
\node at (260,175) {\scalebox{1.2}{$A_3$}};
\node at (330,250) {\scalebox{1.2}{$A_0$}};
\node at (400,120) {\scalebox{1.2}{$A_1$}};

\node at (100,435) {\scalebox{1.2}{$B_0$}};
\node at (240,365) {\scalebox{1.2}{$B_1$}};
\node at (165,295) {\scalebox{1.2}{$B_2$}};
\node at (310,225) {\scalebox{1.2}{$B_3$}};
\node at (230,155) {\scalebox{1.2}{$B_0$}};
\node at (370,90) {\scalebox{1.2}{$B_1$}};

\node at (90,310) {\scalebox{1.2}{$C_1$}};
\node at (160,165) {\scalebox{1.2}{$C_0$}};
\node at (155,380) {\scalebox{1.2}{$C_2$}};
\node at (230,240) {\scalebox{1.2}{$C_3$}};
\node at (305,300) {\scalebox{1.2}{$C_0$}};
\node at (365,175) {\scalebox{1.2}{$C_1$}};
\node at (290,110) {\scalebox{1.2}{$C_2$}};

  \begin{scope}[cm={{0.86169,0.0,0.0,0.86169,(-84.46029,-147.85209)}}]
    \path[draw=black,line join=bevel,line cap=rect,miter limit=4.00,draw
      opacity=0.000,fill opacity=0.627,line width=0.800pt] (580.0000,552.3622) --
      (500.0000,552.3622) -- (500.0000,472.3622) -- (580.0000,393.1662);
    \begin{scope}[shift={(-318.0466,-86.52797)}]
      \path[draw=black,line join=round,line cap=round,even odd rule,line
        width=0.800pt] (490.0000,667.1122) -- (410.0000,667.1122);
    \end{scope}
    \begin{scope}[shift={(-237.84955,-246.47929)}]
      \path[draw=black,line join=round,line cap=round,even odd rule,line
        width=0.800pt] (490.0000,667.1122) -- (410.0000,667.1122);
    \end{scope}
    \begin{scope}[shift={(-316.42646,152.99713)}]
      \path[draw=black,line join=round,line cap=round,even odd rule,line
        width=0.800pt] (410.0000,587.3622) -- (490.0000,507.3622);
    \end{scope}
    \begin{scope}[shift={(-236.22941,-6.95419)}]
      \path[draw=black,line join=round,line cap=round,even odd rule,line
        width=0.800pt] (490.0000,667.1122) -- (410.0000,667.1122) --
        (410.0000,587.3622) -- (490.0000,507.3622);
    \end{scope}
    \begin{scope}[shift={(-156.03236,-166.90552)}]
      \path[draw=black,line join=round,line cap=round,even odd rule,line
        width=0.800pt] (490.0000,667.1122) -- (410.0000,667.1122) --
        (410.0000,587.3622) -- (490.0000,507.3622);
    \end{scope}
    \begin{scope}[shift={(-75.83531,-326.85684)}]
      \path[draw=black,line join=round,line cap=round,even odd rule,line
        width=0.800pt] (490.0000,667.1122) -- (410.0000,667.1122) --
        (410.0000,587.3622);
    \end{scope}
    \begin{scope}[shift={(-154.41222,72.61958)}]
      \path[draw=black,line join=round,line cap=round,even odd rule,line
        width=0.800pt] (410.0000,667.1122) -- (410.0000,587.3622) --
        (490.0000,507.3622);
    \end{scope}
    \begin{scope}[shift={(-74.21516,-87.33174)}]
      \path[draw=black,line join=round,line cap=round,even odd rule,line
        width=0.800pt] (490.0000,667.1122) -- (410.0000,667.1122) --
        (410.0000,587.3622) -- (490.0000,507.3622);
    \end{scope}
    \begin{scope}[shift={(5.98189,-247.28307)}]
      \path[draw=black,line join=round,line cap=round,even odd rule,line
        width=0.800pt] (490.0000,667.1122) -- (410.0000,667.1122) --
        (410.0000,587.3622) -- (490.0000,507.3622);
    \end{scope}
    \begin{scope}[shift={(86.17894,-407.23439)}]
      \path[draw=black,line join=round,line cap=round,even odd rule,line
        width=0.800pt] (490.0000,667.1122) -- (410.0000,667.1122) --
        (410.0000,587.3622);
    \end{scope}
    \begin{scope}[shift={(7.60203,-7.75797)}]
      \path[draw=black,line join=round,line cap=round,even odd rule,line
        width=0.800pt] (410.0000,667.1122) -- (410.0000,587.3622) --
        (490.0000,507.3622);
    \end{scope}
    \begin{scope}[shift={(87.79908,-167.70929)}]
      \path[draw=black,line join=round,line cap=round,even odd rule,line
        width=0.800pt] (490.0000,667.1122) -- (410.0000,667.1122) --
        (410.0000,587.3622) -- (490.0000,507.3622);
    \end{scope}
    \begin{scope}[shift={(167.99613,-327.66062)}]
      \path[draw=black,line join=round,line cap=round,even odd rule,line
        width=0.800pt] (490.0000,667.1122) -- (410.0000,667.1122) --
        (410.0000,587.3622) -- (490.0000,507.3622);
    \end{scope}
  \end{scope}

\end{tikzpicture}

	}
  \caption{$2\times 2$ hexagon momentum polytope of the fundamental domain in $U(\Fd)$.}
	\label{fig:2hexagon}
\end{figure}

We can choose a basis for the homology classes of curves of $\Fd$ given by the 6 curves $A_0, A_1, B_0, B_1, C_0$, and $C_1$, as labelled in Figure~\ref{fig:2hexagon}. This can be shown using simple applications of the standard hexagon relations (Eq.~\ref{eq:hex}) as follows.

First observe that the sum of the the $C$ curves in each row is constant, as is the sum in each column,
\begin{equation}
\begin{aligned}
C_0 + C_2 &= C_1+C_3,\\
C_0 + C_1 &= C_2+C_3. 
\label{eq:hexagon22a}
\end{aligned}
\end{equation}

This follows by combining two hexagon relations to the bottom row of hexagons in Fig.~\ref{fig:2hexagon}. For example, 
\begin{equation*}
\begin{aligned}
A_0+C_1&=A_2+C_2,\\ 
A_0+C_0&=A_2+C_3.
\end{aligned}
\end{equation*}
yields $C_0 + C_2 = C_1+C_3$. The other relation is derived similarly.

We can also write the sum of all the diagonal curves in two ways, grouped as rows or as columns. In this case when $I=J=2$, we have:
\begin{equation*}
\begin{aligned}
(C_0 + C_1) + (C_2 + C_3) = (C_0 + C_2) + (C_1 + C_3),
\end{aligned}
\end{equation*}

Together with the previous Eq.~\ref{eq:hexagon22a}, this implies that
\begin{equation*}
\begin{aligned}
C_1=C_2,\\
C_0=C_3.
\end{aligned}
\end{equation*}

In a similar fashion, it is easy to deduce that 
\begin{equation*}
\begin{aligned}
A_0=A_2, \\
A_1=A_3, \\
B_0=B_2, \\
B_1=B_3.
\end{aligned}
\end{equation*}
We can thus choose a basis for the homology classes of curves given by the 6 curves $A_0, A_1, B_0, B_1, C_0, C_1$.

We are interested in sheaves invariant under the action of the torus $T$, so their support must be contained in the $T$-torus fixed curves. We will assume from now on that the support curve $C$ of our sheaf has $\deg B_0=1$ and $\deg B_1=0$ so that it is in the homology class:
\begin{equation}
[C] =\beta \in \sum_{i=0,1}{a_i[A_i]} +[B_0]+ \sum_{j=0,1}{c_j[C_j]}.
\label{eq:2x2curveclass}
\end{equation}
Recall there is a 1-1 correspondence between the $\pi_1(\Fd)$-equivariant sheaves on $\UFm$ and the $P$-fixed sheaves on $\Fm$ up to deck transformations \cite[Proposition~22]{nina19}. To remove the ambiguity, we will further require that the corresponding equivariant sheaf in $\UFm$ has support curve $\CCC\subset\UFm$, whose reduced irreducible component covering the curve in class $[B_0]$ is in our chosen fundamental domain $D$. In this example, notice that there are two possible choices for $B_0$, since $[B_0]=[B_2]$.

The irreducible components of the curve $\CCC$ may be nonreduced, so we need to keep track of the multiplicity of each component to determine the curve class of $pr(\CCC)$. We let the variable $r_i$ track the number of curves which cover $A_i$, $i\in\{0,1\}$, and $s_j$ track the number of curves which cover $C_j$, $j\in\{0,1\}$. A typical curve $\CCC$ which covers a curve in class Eq.~(\ref{eq:2x2curveclass}) is pictured in Figure~\ref{fig:2hexcurve}.
\begin{figure}[th!]
  \centering
	\scalebox{.8}{

\begin{tikzpicture}[y=0.80pt, x=0.80pt, yscale=-1.000000, xscale=1.000000, inner sep=0pt, outer sep=0pt]
\node at (180,315) {\scalebox{1}{$\mathcal{B}_0$}};

\node at (70,440) {\scalebox{1}{$r_0$}};
\node at (90,400) {\scalebox{1}{$s_1$}};
\node at (110,365) {\scalebox{1}{$r_1$}};
\node at (130,320) {\scalebox{1}{$s_0$}};

\node at (150,285) {\scalebox{1}{$r_0$}};
\node at (170,240) {\scalebox{1}{$s_1$}};
\node at (190,205) {\scalebox{1}{$r_1$}};
\node at (210,160) {\scalebox{1}{$s_0$}};

\node at (230,125) {\scalebox{1}{$r_0$}};
\node at (250,85) {\scalebox{1}{$s_1$}};
\node at (270,45) {\scalebox{1}{$r_1$}};

\node at (150,445) {\scalebox{1}{$s_1$}};
\node at (170,400) {\scalebox{1}{$r_0$}};
\node at (185,370) {\scalebox{1}{$s_0$}};
\node at (210,325) {\scalebox{1}{$r_1$}};

\node at (225,290) {\scalebox{1}{$s_1$}};
\node at (250,250) {\scalebox{1}{$r_0$}};
\node at (265,215) {\scalebox{1}{$s_0$}};
\node at (285,165) {\scalebox{1}{$r_1$}};

\node at (305,135) {\scalebox{1}{$s_1$}};
\node at (325,85) {\scalebox{1}{$r_0$}};
\node at (340,55) {\scalebox{1}{$s_0$}};

  \begin{scope}[shift={(177.98554,-588.65545)}]
  \end{scope}
  \path[draw=black,line join=round,line cap=round,even odd rule,line
    width=0.393pt] (40.7352,500.9359) -- (79.9947,461.6764);
  \begin{scope}[cm={{0.49074,0.0,0.0,0.49074,(-160.46962,212.69177)}}]
    \path[draw=black,line join=round,line cap=round,even odd rule,line
      width=0.800pt] (410.0000,587.3622) -- (490.0000,507.3622);
  \end{scope}
  \begin{scope}[cm={{0.49074,0.0,0.0,0.49074,(-121.11344,134.19671)}}]
    \path[draw=black,line join=round,line cap=round,even odd rule,line
      width=0.800pt] (410.0000,667.1122) -- (410.0000,587.3622) --
      (490.0000,507.3622);
  \end{scope}
  \begin{scope}[cm={{0.49074,0.0,0.0,0.49074,(-81.75727,55.70164)}}]
    \path[draw=black,line join=round,line cap=round,even odd rule,line
      width=0.800pt] (410.0000,667.1122) -- (410.0000,587.3622) --
      (490.0000,507.3622);
  \end{scope}
  \begin{scope}[cm={{0.49074,0.0,0.0,0.49074,(-42.40109,-22.79342)}}]
    \path[draw=black,line join=round,line cap=round,even odd rule,line
      width=0.800pt] (490.0000,667.1122) -- (410.0000,667.1122) --
      (410.0000,587.3622) -- (490.0000,507.3622);
  \end{scope}
  \begin{scope}[cm={{0.49074,0.0,0.0,0.49074,(-3.04491,-101.28849)}}]
    \path[draw=black,line join=round,line cap=round,even odd rule,line
      width=0.800pt] (410.0000,667.1122) -- (410.0000,587.3622) --
      (490.0000,507.3622);
  \end{scope}
  \begin{scope}[cm={{0.49074,0.0,0.0,0.49074,(36.31126,-179.78356)}}]
    \path[draw=black,line join=round,line cap=round,even odd rule,line
      width=0.800pt] (410.0000,667.1122) -- (410.0000,587.3622) --
      (490.0000,507.3622);
  \end{scope}
  \begin{scope}[cm={{0.49074,0.0,0.0,0.49074,(75.66744,-258.27862)}}]
    \path[draw=black,line join=round,line cap=round,even odd rule,line
      width=0.800pt] (410.0000,667.1122) -- (410.0000,587.3622);
  \end{scope}
  \begin{scope}[cm={{0.49074,0.0,0.0,0.49074,(-80.96219,173.24701)}}]
    \path[draw=black,line join=round,line cap=round,even odd rule,line
      width=0.800pt] (410.0000,667.1122) -- (410.0000,587.3622) --
      (490.0000,507.3622);
  \end{scope}
  \begin{scope}[cm={{0.49074,0.0,0.0,0.49074,(-41.60602,94.75195)}}]
    \path[draw=black,line join=round,line cap=round,even odd rule,line
      width=0.800pt] (410.0000,667.1122) -- (410.0000,587.3622) --
      (490.0000,507.3622);
  \end{scope}
  \begin{scope}[cm={{0.49074,0.0,0.0,0.49074,(-2.24984,16.25688)}}]
    \path[draw=black,line join=round,line cap=round,even odd rule,line
      width=0.800pt] (410.0000,667.1122) -- (410.0000,587.3622) --
      (490.0000,507.3622);
  \end{scope}
  \begin{scope}[cm={{0.49074,0.0,0.0,0.49074,(37.10634,-62.23818)}}]
    \path[draw=black,line join=round,line cap=round,even odd rule,line
      width=0.800pt] (410.0000,667.1122) -- (410.0000,587.3622) --
      (490.0000,507.3622);
  \end{scope}
  \begin{scope}[cm={{0.49074,0.0,0.0,0.49074,(76.46251,-140.73325)}}]
    \path[draw=black,line join=round,line cap=round,even odd rule,line
      width=0.800pt] (410.0000,667.1122) -- (410.0000,587.3622) --
      (490.0000,507.3622);
  \end{scope}
  \begin{scope}[cm={{0.49074,0.0,0.0,0.49074,(115.81869,-219.22831)}}]
    \path[draw=black,line join=round,line cap=round,even odd rule,line
      width=0.800pt] (410.0000,667.1122) -- (410.0000,587.3622) --
      (490.0000,507.3622);
  \end{scope}
  \path[draw=black,fill=black,line join=round,line cap=round,miter limit=4.00,draw
    opacity=0.000,line width=0.592pt]
    (310.8821,20.1203)arc(-0.000:89.774:1.608)arc(89.773:179.547:1.608)arc(179.547:269.321:1.608)arc(269.321:359.094:1.608);
  \path[draw=black,fill=black,line join=round,line cap=round,miter limit=4.00,draw
    opacity=0.000,line width=0.592pt]
    (310.8821,28.2830)arc(0.000:89.774:1.608)arc(89.773:179.547:1.608)arc(179.547:269.321:1.608)arc(269.321:359.094:1.608);
  \path[draw=black,fill=black,line join=round,line cap=round,miter limit=4.00,draw
    opacity=0.000,line width=0.592pt]
    (310.8821,36.4458)arc(0.000:89.774:1.608)arc(89.773:179.547:1.608)arc(179.547:269.321:1.608)arc(269.321:359.094:1.608);
  \path[draw=black,fill=black,line join=round,line cap=round,miter limit=4.00,draw
    opacity=0.000,line width=0.592pt]
    (89.8525,488.9635)arc(-0.000:89.773:1.608)arc(89.773:179.547:1.608)arc(179.547:269.321:1.608)arc(269.320:359.094:1.608);
  \path[draw=black,fill=black,line join=round,line cap=round,miter limit=4.00,draw
    opacity=0.000,line width=0.592pt]
    (89.8525,497.1263)arc(-0.000:89.773:1.608)arc(89.773:179.547:1.608)arc(179.547:269.321:1.608)arc(269.320:359.094:1.608);
  \path[draw=black,fill=black,line join=round,line cap=round,miter limit=4.00,draw
    opacity=0.000,line width=0.592pt]
    (89.8525,505.2890)arc(-0.000:89.773:1.608)arc(89.773:179.547:1.608)arc(179.547:269.321:1.608)arc(269.320:359.094:1.608);

\end{tikzpicture}
	
	}
  \caption{A typical torus fixed curve in $U(\Fd)$ which covers a curve with $\deg [B]=\deg [B_0]=1$. Here $r_i$ and $s_j$ are used to track the multiplicity of curves which cover $A_i$ and $C_j$ curves, respectively.}
	\label{fig:2hexcurve}
\end{figure}

\subsection{Translating to combinatorics}
We explain the details of our method of converting the count of invariant stable sheaves into a combinatorics problem.

Recall from the discussion in Section~\ref{subsec:method} the naive Euler characteristic $\widetilde{n}^0_{\beta}(\Ld) = e(M^{\Ld}_{\beta})$ for curves in class $\beta$ of the form Eq.~(\ref{eq:2x2curveclass}) equals a count of $T$-torus invariant structure sheaves of genus 0 curves on the universal cover $\CCC\subset\UFm$ that cover class $\beta$ \cite[Proposition 23]{nina19}, subject to the certain conditions  \cite[Proposition 31]{nina19} on their multiplicity that we explain below.

First we introduce some terminology. We will refer to irreducible components of $\CCC$ as edges, and the intersection of two or more edges as vertices. Let $\mathcal{B}_0\subset\CCC$ be the edge that covers class $[B_0]$ and lies in the fundamental domain $D$ by assumption.

We will call the union of $\mathcal{B}_0$ with any one of the four disjoint subcurves of $\overline{\CCC\backslash \mathcal{B}_0}$ a branch of $\CCC$. Then the curve counts can be done on each branch separately, and will be the same on each branch, up to relabeling.

The edge $\mathcal{B}_0$ is the intersection of two irreducible surface component hexagons isomorphic to $\Xi$ in $\UFm$. Let $S$ be one of these and let $g$ be the deck transformation that translates $S$ into the other. The hexagons $g^m S$, $m\in\ZZ$, in the orbit of $S$ under the group of deck transformations $\langle g \rangle\cong\ZZ$ will be called inside hexagons. Any other hexagons will be called outside hexagons. 

Any edge of $\overline{\CCC\backslash \mathcal{B}_0}$ covers $A_i$ or $C_j$, $i,j\in\{0,1\}$, and will be the intersection of an inside hexagon and an outside one. These $T$-invariant edges can can have monomial thickening in these two directions. Any thickenings in the direction of the inside hexagon will be called inside thickenings, and those in the direction of the outside hexagon will be called outside thickenings. However, because of stability and the requirement of the Euler characteristic of $\OO_{[\CCC]}$ to be 1, the possible thickenings can only be of a particular form. 

Thickenings on the edges that cover $A_i$ or $C_j$ are subject to the following properties \cite[Proposition 31]{nina19}:
\begin{enumerate}
	\item Inside thickenings of any edge that intersects $\mathcal{B}_0$ is unrestricted.
	\item All nonzero outside thickenings must be 1.
	\item Inside thickenings are non-increasing on components along a branch in the direction moving away from $\mathcal{B}_0$.
	\item Inside thickenings for two adjacent edges contained in a common inside hexagon can either be the same or differ by one.  
\end{enumerate} 
We can interpret the inside multiplicity of each edge as length of a part in a partition. These constraints are independent on each branch, so we examine one branch at a time. Along each branch, the non-increasing length condition says that the allowed multiplicities of edges form a Young diagram. If we examine the conjugate partition of the branch, the fourth condition can be interpreted as saying that any odd parts that appear in the conjugate partition are distinct, with no restriction on the even parts. 

The generating function that counts the number of partitions $p(n)$ with odd parts distinct can be written as the product of generating functions for partitions with arbitrary even parts with that of partitions with unique odd parts: 
\[
\sum{p(n)q^n}=\prod{\frac{1}{1-q^{2n}}}\prod{(1 +q^{2n-1})}
\] 
We must further refine the parts, because we have four curve classes $A_i$ and $C_j$, for $i,j\in\{0,1\}$ to keep track of. In other words, we need to keep track of the residue classes mod 4 in our partition.  

\begin{figure}[th!]
  \centering
	\scalebox{.8}{

\begin{tikzpicture}[y=0.80pt, x=0.80pt, yscale=-1.000000, xscale=1.000000, inner sep=0pt, outer sep=0pt]

\node at (180,315) {\scalebox{1}{$\mathcal{B}_0$}};

\node at (225,290) {\scalebox{1}{$s_1$}};
\node at (250,250) {\scalebox{1}{$r_0$}};
\node at (265,215) {\scalebox{1}{$s_0$}};
\node at (285,165) {\scalebox{1}{$r_1$}};

\node at (305,135) {\scalebox{1}{$s_1$}};
\node at (325,85) {\scalebox{1}{$r_0$}};
\node at (340,55) {\scalebox{1}{$s_0$}};

  \path[draw=black,line join=round,line cap=round,even odd rule,line
    width=0.393pt] (198.0632,304.5875) -- (158.8037,304.5875);
  \path[draw=black,line join=round,line cap=round,even odd rule,line
    width=0.393pt] (198.9550,304.5010) -- (238.2145,265.2416);
  \path[draw=black,line join=round,line cap=round,even odd rule,line
    width=0.393pt] (238.3112,265.1428) -- (238.3112,226.0060) --
    (277.5706,186.7465);
  \path[draw=black,line join=round,line cap=round,even odd rule,line
    width=0.393pt] (277.6673,186.6477) -- (277.6673,147.5109) --
    (316.9268,108.2514);
  \path[draw=black,line join=round,line cap=round,even odd rule,line
    width=0.393pt] (317.0235,108.1526) -- (317.0235,69.0158) --
    (356.2830,29.7564);
  \path[draw=black,fill=black,line join=round,line cap=round,miter limit=4.00,draw
    opacity=0.000,line width=0.592pt]
    (310.8821,20.1203)arc(-0.000:89.774:1.608)arc(89.773:179.547:1.608)arc(179.547:269.321:1.608)arc(269.321:359.094:1.608);
  \path[draw=black,fill=black,line join=round,line cap=round,miter limit=4.00,draw
    opacity=0.000,line width=0.592pt]
    (310.8821,28.2830)arc(0.000:89.774:1.608)arc(89.773:179.547:1.608)arc(179.547:269.321:1.608)arc(269.321:359.094:1.608);
  \path[draw=black,fill=black,line join=round,line cap=round,miter limit=4.00,draw
    opacity=0.000,line width=0.592pt]
    (310.8821,36.4458)arc(0.000:89.774:1.608)arc(89.773:179.547:1.608)arc(179.547:269.321:1.608)arc(269.321:359.094:1.608);

\end{tikzpicture}

	}
  \caption{Detail of the northeast branch of the curve shown in Figure~\ref{fig:2hexcurve}.}
	\label{fig:2hexcurveNE}
\end{figure}
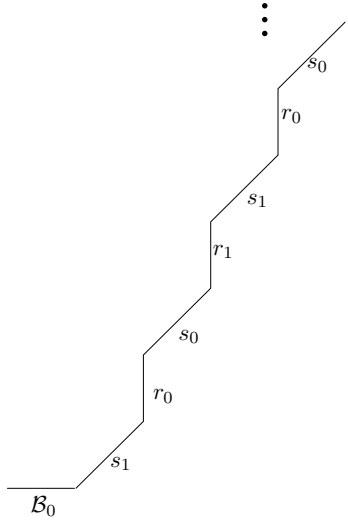

Consider for example the northeast branch of the curve shown above, which we reproduce in Figure~\ref{fig:2hexcurveNE}. Suppose we number the edges consecutively, starting with the first edge $e_1$ that intersects $\mathcal{B}_0$. Then the every odd-numbered edge will contribute to an odd part, and the even-numbered ones to an even part. The first edge $e_1$ in this example curve covers $C_1$, and we assign the variable $s_1$ to track this curve. The second edge $e_2$ covers $A_0$ and we use the variable $r_0$ to track this.

So we refine the generating function above, and replace powers of the variable $q$, by
\begin{align*}
q^1 &\mapsto s_1 \\
q^2 &\mapsto s_1 r_0\\
q^3 &\mapsto s_1 r_0 s_0 \\
q^4 &\mapsto s_1 r_0 s_0 r_1,
\end{align*}
and for higher powers of $q$, we continue the pattern, so that
\begin{align*}
q^{4m+i}=q^{4m}q^i &\mapsto (s_1 r_0 s_0 r_1)^m q^i
\end{align*}

We also use the notation 
\[
Q\coloneqq r_0s_0r_1s_1.
\]

Now, the generating function that counts the number of partitions with odd parts distinct can be expressed for the northeast branch as:
\begin{align*}
&\frac{(1+s_1)(1+s_0r_0s_1)(1+Qs_1)(1+Qs_0r_0s_1)\cdots}{(1-r_0s_1)(1-Q)(1-Qr_0s_1)(1-Q^2)\cdots} \\
&=\frac{ (1+s_1)\prod\limits_{m=1}^{\infty}{(1+Q^ms_1)(1+Q^mr_1^{-1})} }{(1-r_0s_1)\prod\limits_{m=1}^{\infty}{(1-Q^m)(1-r_0s_1 Q^m)} }
\end{align*}

We do this for each branch and multiply the contribution from all four branches. Notice also that there are two distinct possible locations for the curve $B_0$. However, they give the same contribution to the generating function, since the four branches in either location consist of the same sequence of curves, up to renaming of the branches.

Hence, using the identities (\ref{eq:Jphi-identities}), the partition function can be expressed in terms of the theta function $\thetaQ{p}$ as:
\begin{align}
&2i\etaQ^{-6}\,\frac{\thetaQ{-r_0}\thetaQ{-s_0}\thetaQ{-r_1}\thetaQ{-s_1} } { \thetaQ{r_0s_0}\thetaQ{r_1s_1}}.\nonumber \\ 
\end{align}
or in terms of the weak Jacobi form $\JphiQ{p}$ as:
\begin{align}
&2\left\{\frac{\JphiQ{-r_0}\JphiQ{-s_0}\JphiQ{-r_1}\JphiQ{-s_1} } { \JphiQ{r_0s_0}\JphiQ{r_1s_1} } \right\}^{1/2} \nonumber \\ 
\end{align}
As we explained in the introduction, this count of the fixed points corresponds to the naive Euler characteristic contribution, $\naivel$, 
\begin{align}
\sum_{\substack{a_0, a_1\\c_0, c_1}} {\widetilde{n}^0_{\beta_{(\textbf{a},\textbf{c})}}(\Lm) r_0^{a_0}r_1^{a_1} s_0^{c_0} s_1^{c_1}} =
2\left\{\frac{\JphiQ{-r_0}\JphiQ{-s_0}\JphiQ{-r_1}\JphiQ{-s_1} } { \JphiQ{r_0s_0}\JphiQ{r_1s_1} } \right\}^{1/2}.\nonumber \\ 
\end{align}

However, the Behrend function weighting amounts to a sign that depends on the degree of the curve class \cite[Remark 33]{nina19}:
\[
\widetilde{n}^0_{\beta_{(\textbf{a},\textbf{c})}}(\Lm) = (-1)^{a_0 + a_1 + c_0 + c_1}{n}^0_{\beta_{(\textbf{a},\textbf{c})}}(\Lm).
\]
We can incorporate this sign by replacing our tracking variables by their negatives.

Hence, we have shown the following.
\begin{theorem}\label{thm:22}
Fix a curve class $\beta_{(\textbf{a},\,\textbf{c})}$ in the local multi-Banana $\Lm=\Lvw{2}{2}$,
\begin{gather*}
\beta_{(\textbf{a},\,\textbf{c})} = {a_0}[A_0] + {a_1}[A_1] +{c_0}[C_0] + {c_1}[C_1] + [B_0],\\
\textbf{ a} = (a_0,a_1), \textbf{ c} =(c_0,c_1) \in \ZZ_{\geq 0}^2.
\end{gather*}
Then the genus 0 Gopakumar-Vafa invariants $n^0_{\beta_{(\textbf{a},\textbf{c})}}(\Lm)$ are given by the following:
\[
\sum_{\substack{a_0, a_1\\c_0, c_1}}{n^0_{\beta_{(\textbf{a},\textbf{c})}}(\Lm) r_0^{a_0}r_1^{a_1} s_0^{c_0} s_1^{c_1}} =
2\left\{\frac{\JphiQ{r_0}\JphiQ{s_0}\JphiQ{r_1}\JphiQ{s_1} } { \JphiQ{r_0s_0}\JphiQ{r_1s_1} } \right\}^{1/2},
\]
where we use the notation
\begin{align*}
Q&\coloneqq r_0 r_1 s_0 s_1\\
\JphiQ{p}&\coloneqq \Jphi{Q}{p}.
\end{align*}
and $\JphiQ{p}=\Jphi{Q}{p}$ is the unique weak Jacobi form of weight -2 and index 1:
\[
\JphiQ{p} = p^{-1}(1-p)^2\prod_{m=1}^{\infty}{\frac{(1-Q^{m}p^{-1})^2 (1-Q^{m}p)^2}{(1-Q^m)^4}},
\]
\[Q=\exp(2\pi i\tau),  \qquad p = \exp(2\pi iz), \qquad  (\tau,z)\in \HH\times\CC.\]

\end{theorem}

%%%%%%%%%%%%%%%%%%%%%%%%%%%%%%%%%%%%%%%%%%%
%\input{./section4}
\section{Case $1\times w$}\label{sec:1w}
In this section, we look in detail at the case of $\Fvw{1}{w}$, when $v=1$ and $w\geq 1$.

In this case, the fundamental domain in $U(\Fvw{1}{w})$ has a $\ZZ\times w\ZZ $ periodicity in the universal cover. Its momentum polytope is given by $1\times w$ hexagons. Using the hexagon relations (Eq.~\ref{eq:hex}) as in the previous section, it is easy to see there is only one horizontal curve class, which we call $B$, and one diagonal curve class, which we call $C$. There are $w$ distinct vertical curve classes, $A_i, 0\leq i \leq w-1$, (Figure~\ref{fig:1xwhexagon}).

\begin{figure}[th!]
  \centering
	\scalebox{.5}{

\begin{tikzpicture}[y=0.80pt, x=0.80pt, yscale=-1.000000, xscale=1.000000, inner sep=0pt, outer sep=0pt]

\node at (70,590) {\scalebox{1.2}{$A_0$}};
\node at (120,480) {\scalebox{1.2}{$A_1$}};
\node at (175,375) {\scalebox{1.2}{$A_2$}};
\node at (272,160) {\scalebox{1.2}{$A_{w-1}$}};

\node at (200,535) {\scalebox{1.2}{$A_0$}};
\node at (255,425) {\scalebox{1.2}{$A_1$}};
\node at (305,320) {\scalebox{1.2}{$A_2$}};
\node at (423,105) {\scalebox{1.2}{$A_{w-1}$}};

\node at (105,625) {\scalebox{1.2}{$B$}};
\node at (155,520) {\scalebox{1.2}{$B$}};
\node at (215,395) {\scalebox{1.2}{$B$}};
\node at (325,175) {\scalebox{1.2}{$B$}};
\node at (375,75) {\scalebox{1.2}{$B$}};

\node at (100,525) {\scalebox{1.2}{$C$}};
\node at (155,420) {\scalebox{1.2}{$C$}};
\node at (200,320) {\scalebox{1.2}{$C$}};
\node at (310,105) {\scalebox{1.2}{$C$}};

\node at (170,595) {\scalebox{1.2}{$C$}};
\node at (220,490) {\scalebox{1.2}{$C$}};
\node at (285,375) {\scalebox{1.2}{$C$}};
\node at (385,165) {\scalebox{1.2}{$C$}};

  \begin{scope}[shift={(-494.51803,320.69501)}]
  \end{scope}
  \begin{scope}[cm={{0.66421,0.0,0.0,0.66421,(-245.89801,118.92373)}}]
    \path[draw=black,line join=round,line cap=round,even odd rule,line
      width=0.800pt] (490.0000,667.1122) -- (410.0000,667.1122);
  \end{scope}
  \begin{scope}[cm={{0.66421,0.0,0.0,0.66421,(-192.63002,12.68185)}}]
    \path[draw=black,line join=round,line cap=round,even odd rule,line
      width=0.800pt] (490.0000,667.1122) -- (410.0000,667.1122);
  \end{scope}
  \begin{scope}[cm={{0.66421,0.0,0.0,0.66421,(-139.36202,-93.56005)}}]
    \path[draw=black,line join=round,line cap=round,even odd rule,line
      width=0.800pt] (490.0000,667.1122) -- (410.0000,667.1122);
  \end{scope}
  \begin{scope}[cm={{0.66421,0.0,0.0,0.66421,(-86.09403,-199.80194)}}]
  \end{scope}
  \begin{scope}[cm={{0.66421,0.0,0.0,0.66421,(-32.82604,-306.04381)}}]
    \path[draw=black,line join=round,line cap=round,even odd rule,line
      width=0.800pt] (490.0000,667.1122) -- (410.0000,667.1122);
  \end{scope}
  \begin{scope}[cm={{0.66421,0.0,0.0,0.66421,(-244.8219,278.01961)}}]
    \path[draw=black,line join=round,line cap=round,even odd rule,line
      width=0.800pt] (410.0000,587.3622) -- (490.0000,507.3622);
  \end{scope}
  \begin{scope}[cm={{0.66421,0.0,0.0,0.66421,(-191.5539,171.77774)}}]
    \path[draw=black,line join=round,line cap=round,even odd rule,line
      width=0.800pt] (490.0000,667.1122) -- (410.0000,667.1122) --
      (410.0000,587.3622) -- (490.0000,507.3622);
  \end{scope}
  \begin{scope}[cm={{0.66421,0.0,0.0,0.66421,(-138.28591,65.53585)}}]
    \path[draw=black,line join=round,line cap=round,even odd rule,line
      width=0.800pt] (490.0000,667.1122) -- (410.0000,667.1122) --
      (410.0000,587.3622) -- (490.0000,507.3622);
  \end{scope}
  \begin{scope}[cm={{0.66421,0.0,0.0,0.66421,(-85.01792,-40.70604)}}]
    \path[draw=black,line join=round,line cap=round,even odd rule,line
      width=0.800pt] (490.0000,667.1122) -- (410.0000,667.1122) --
      (410.0000,587.3622) -- (490.0000,507.3622);
  \end{scope}
  \begin{scope}[cm={{0.66421,0.0,0.0,0.66421,(-31.74992,-146.94792)}}]
    \path[draw=black,line join=round,line cap=round,even odd rule,line
      width=0.800pt] (410.0000,587.3622) -- (490.0000,507.3622);
  \end{scope}
  \begin{scope}[cm={{0.66421,0.0,0.0,0.66421,(21.51807,-253.1898)}}]
    \path[draw=black,line join=round,line cap=round,even odd rule,line
      width=0.800pt] (490.0000,667.1122) -- (410.0000,667.1122) --
      (410.0000,587.3622) -- (490.0000,507.3622);
  \end{scope}
  \begin{scope}[cm={{0.66421,0.0,0.0,0.66421,(74.78607,-359.43171)}}]
    \path[draw=black,line join=round,line cap=round,even odd rule,line
      width=0.800pt] (490.0000,667.1122) -- (410.0000,667.1122) --
      (410.0000,587.3622);
  \end{scope}
  \begin{scope}[cm={{0.66421,0.0,0.0,0.66421,(-137.20979,224.63174)}}]
    \path[draw=black,line join=round,line cap=round,even odd rule,line
      width=0.800pt] (410.0000,667.1122) -- (410.0000,587.3622) --
      (490.0000,507.3622);
  \end{scope}
  \begin{scope}[cm={{0.66421,0.0,0.0,0.66421,(-83.94179,118.38986)}}]
    \path[draw=black,line join=round,line cap=round,even odd rule,line
      width=0.800pt] (490.0000,667.1122) -- (410.0000,667.1122) --
      (410.0000,587.3622) -- (490.0000,507.3622);
  \end{scope}
  \begin{scope}[cm={{0.66421,0.0,0.0,0.66421,(-30.67379,12.14796)}}]
    \path[draw=black,line join=round,line cap=round,even odd rule,line
      width=0.800pt] (490.0000,667.1122) -- (410.0000,667.1122) --
      (410.0000,587.3622) -- (490.0000,507.3622);
  \end{scope}
  \begin{scope}[cm={{0.66421,0.0,0.0,0.66421,(22.5942,-94.09393)}}]
    \path[draw=black,line join=round,line cap=round,even odd rule,line
      width=0.800pt] (490.0000,667.1122) -- (410.0000,667.1122) --
      (410.0000,587.3622);
  \end{scope}
  \begin{scope}[cm={{0.66421,0.0,0.0,0.66421,(75.86219,-200.33584)}}]
    \path[draw=black,line join=round,line cap=round,even odd rule,line
      width=0.800pt] (410.0000,667.1122) -- (410.0000,587.3622) --
      (490.0000,507.3622);
  \end{scope}
  \begin{scope}[cm={{0.66421,0.0,0.0,0.66421,(129.13018,-306.57771)}}]
    \path[draw=black,line join=round,line cap=round,even odd rule,line
      width=0.800pt] (490.0000,667.1122) -- (410.0000,667.1122) --
      (410.0000,587.3622) -- (490.0000,507.3622);
  \end{scope}
  \path[draw=black,fill=black,line join=round,line cap=round,miter limit=4.00,draw
    opacity=0.000,line width=0.917pt]
    (284.2347,255.3929)arc(0.000:89.773:2.491)arc(89.773:179.547:2.491)arc(179.547:269.321:2.491)arc(269.321:359.094:2.491);
  \path[draw=black,fill=black,line join=round,line cap=round,miter limit=4.00,draw
    opacity=0.000,line width=0.917pt]
    (284.2347,268.0411)arc(0.000:89.773:2.491)arc(89.773:179.547:2.491)arc(179.547:269.321:2.491)arc(269.321:359.094:2.491);
  \path[draw=black,fill=black,line join=round,line cap=round,miter limit=4.00,draw
    opacity=0.000,line width=0.917pt]
    (284.2347,280.6894)arc(-0.000:89.773:2.491)arc(89.773:179.547:2.491)arc(179.547:269.321:2.491)arc(269.321:359.094:2.491);

\end{tikzpicture}
	
	}
  \caption{$1\times w$ hexagon momentum polytope of the fundamental domain in $U(\Fvw{1}{w})$.}
	\label{fig:1xwhexagon}
\end{figure}

Let us assume that the support curve $\CCC$ of our sheaf has $\deg B=1$ so that
\[
[\CCC] \in \sum_{i=0}^w{a_i[A_i]} +[B]+ c[C].
\]
 Let $r_i$ track the number of $A_i$ curves and $s$ track the number of $C$ curves.

For this case, we define the variable $Q$ to be: 
\begin{equation}
Q\coloneqq \prod\limits_{i=0}^{w-1}{(r_is)}.
\label{eq:Q1w}
\end{equation}
We will also use the following multi product notation:
\begin{equation}
 R_{a;b}\coloneqq r_a\cdot r_{a+1}\cdot r_{a+2} \cdots r_b \cdot s^{b-a+1}, \quad a\leq b,
\label{eq:multiRab}
\end{equation}
where the subscript of $r$ is interpreted mod $w$:
\[r_{k+w} \coloneqq r_{[k]}, \quad [k]\in \ZZ/w\ZZ.\] 
In particular, 
\[
R_{0;b}\coloneqq \prod\limits_{i=0}^{b}{(r_is)}.
\]
First, suppose the single $B$ curve is located connected to an $A_0$ curve. Then, we can count the number of partitions with odd parts distinct in the same way as before, and the generating function for these configurations is expressed as follows:
\begin{equation}
\begin{aligned}
&(1+s)^2 \prod\limits_{m=1}^{\infty}{\frac{(1+sQ^m)^2(1+s^{-1}Q^m)^2  }{(1-Q^m)^4} }\times 
 \prod\limits_{k=0}^{w-2}\left\{ \frac{(1+sR_{0;k})(1+s^{-1}R_{0;k})}{(1-R_{0;k})^2}\right. \nonumber \\
& \times \left.   \prod\limits_{m=1}^{\infty}{ 
\frac{(1+sR_{0;k}Q^m)(1+s^{-1}R_{0;k}Q^m)(1+sR^{-1}_{0;k}Q^m)(1+s^{-1}R^{-1}_{0;k}Q^m)}{(1-R_{0;k}Q^m)^2(1-R^{-1}_{0;k}Q^m)^2}  } \right\} 
\end{aligned}
\end{equation} 

We can write this more succinctly using the weak Jacobi form $\JphiQ{p}$ as:
\begin{equation}
 (-s)\JphiQ{-s} \prod\limits_{k=0}^{w-2}  
 \left\{\frac{\sqrt {\JphiQ{-sR_{0;k}}}\sqrt {\JphiQ{-sR_{0;k}}}} 
{ \JphiQ{R_{0;k}}} \right\};
\end{equation}

or alternatively, in terms of the theta function $\thetaQ{p}$ function as:
\begin{equation}
 (-s)\JphiQ{-s} \prod\limits_{k=0}^{w-2}  
 \left\{ \frac{\thetaQ{-sR_{0;k}}\thetaQ{-sR^{-1}_{0;k}}} 
 {\thetaQ{R_{0;k}}\thetaQ{R^{-1}_{0;k}}} \right\} .
\end{equation}

Notice that, from Eq.~(\ref{eq:EllipticGenusC2}), the product can be expressed in terms of the equivariant elliptic genus of $\CC^2$,
\begin{equation}
(-s)\JphiQ{-s} \prod\limits_{k=0}^{w-2}  \mathrm{Ell}_{Q,-s}(\CC^2,R_{0;k}) .
\end{equation}

There are $w$ different locations possible for the $B$ curve in the $1\times w$ hexagon, characterized by the choice of which $A_i$ curve, $0\leq i\leq w-1$, that the $B$ curve is connected to. Although the generating function for the partitions with distinct odd parts associated to these other configurations depends on the particular location of $B$, it is easy to see that it differs from the previous formula only by a cyclic shift of indices in the $R_{0;k}$ variable.

The total partition function counts the contribution from all possible locations of the $B$ curve, and is thus expressed as a sum over the generating functions from each location,

\begin{equation}
(-s)\JphiQ{-s} \sum_{i=0}^{w-1} \quad {\prod\limits_{k=i}^{i+w-2}  \mathrm{Ell}_{Q,-s}(\CC^2,R_{i;k}) }.
\label{eq:1xwpartition}
\end{equation} 

As in the previous section, this count of fixed points corresponds to the naive Euler characteristic. To take account of the Behrend function weighting, we can incorporate a sign based on the degree of the curve class by simply replacing our tracking variables by their negatives.

Thus, for the case of $\Lvw{1}{w}$, we have the following partition function. 
%where the notation is defined in (Eq.~\ref{eq:Q1w}) and (Eq.~\ref{eq:multiRab}).

\begin{theorem}\label{thm:1W}
Fix a curve class $\beta_{(\textbf{a},\,c)}$ in the local multi-Banana $\Lm=\Lvw{1}{w}$: 
\begin{gather*}
\beta_{(\textbf{a},\,c)} = \sum\limits_{i=0}^{w-1}{{a_i}[A_i]} +{c}[C] + [B],\\
\textbf{a} = (a_0,\ldots,a_{w-1}) \in \ZZ_{\geq 0}^w,\, c  \in \ZZ_{\geq 0}.
\end{gather*}

Then the genus 0 Gopakumar-Vafa invariants $n^0_{\beta_{(\textbf{a},c)}}(\Lvw{1}{w})$ can be expressed as:
\[
\sum_{\textbf{a}, c} {n^0_{\beta_{(\textbf{a},\textbf{c})}}(\Lm) \textbf{r}^{\textbf{a}}s^{c}} =
s\cdot\JphiQ{s} \sum_{i=0}^{w-1} \quad {\prod\limits_{k=i}^{i+w-2}  \mathrm{Ell}_{Q,s}(\CC^2,R_{i;k}) }.
\]
where $\mathrm{Ell}_{q,y}(\CC^2,t)$ is the equivariant elliptic genus of $\CC^2$, and we use the notation:
\begin{gather*}
\textbf{r}^{\textbf{a}}  \coloneqq r_0^{a_0}r_1^{a_1}\ldots r_{w-1}^{a_{w-1}},\\
Q\coloneqq \prod\limits_{i=0}^{w-1}{(r_is)},\\
R_{a;b}\coloneqq r_a\cdot r_{a+1}\cdot r_{a+2} \cdots r_b \cdot s^{b-a+1}, \quad a\leq b,\\
r_{k+w} \coloneqq r_{[k]}, \quad [k]\in \ZZ/w\ZZ.
\end{gather*} 
\end{theorem}

We also mention that it is possible to choose to fix the degree of the $A$ family of curve classes or the $C$ class to be 1 instead of the $B$ curve. However, in this $1\times w$ case, doing so reduces to the ordinary Banana configuration $\Fb$ case and yields the same formula as in the earlier paper \cite{nina19}.

%%%%%%%%%%%%%%%%%%%%%%%%%%%%%%%%%%%%%%%%%%%

%\begin{comment}
\bibliography{mainbib}

\begin{thebibliography}{10}

\bibitem{behrend09}
Kai Behrend.
\newblock Donaldson-{T}homas type invariants via microlocal geometry.
\newblock {\em Ann. of Math. (2)}, 170(3):1307--1338, 2009.

\bibitem{behrend-fantechi08}
Kai Behrend and Barbara Fantechi.
\newblock Symmetric obstruction theories and {H}ilbert schemes of points on
  threefolds.
\newblock {\em Algebra Number Theory}, 2(3):313--345, 2008.

\bibitem{bryan19}
Jim Bryan.
\newblock The {D}onaldson-{T}homas partition function of the banana manifold.
\newblock {\em Algebr. Geom.}, 8(2):133--170, 2021.

\bibitem{iqbal13}
Stefan Hohenegger and Amer Iqbal.
\newblock {M-strings, elliptic genera and $\mathcal{N} = 4$ string amplitudes}.
\newblock {\em Fortsch. Phys.}, 62:155--206, 2014.

\bibitem{kanazawa-lau19}
Atsushi Kanazawa and Siu-Cheong Lau.
\newblock Local {C}alabi-{Y}au manifolds of type {$\tilde{A}$} via {SYZ} mirror
  symmetry.
\newblock {\em J. Geom. Phys.}, 139:103--138, 2019.

\bibitem{katz08}
Sheldon Katz.
\newblock Genus zero {G}opakumar-{V}afa invariants of contractible curves.
\newblock {\em J. Differential Geom.}, 79(2):185--195, 2008.

\bibitem{katzklemmvafa97}
Sheldon Katz, Albrecht Klemm, and Cumrun Vafa.
\newblock Geometric engineering of quantum field theories.
\newblock {\em Nuclear Phys. B}, 497(1-2):173--195, 1997.

\bibitem{nina19}
Nina {Morishige}.
\newblock Genus zero {G}opakumar-{V}afa invariants of the {B}anana manifold.
\newblock {\em arXiv e-prints}, page arXiv:1909.01540, September 2019.

\bibitem{simpson94}
Carlos~T. Simpson.
\newblock Moduli of representations of the fundamental group of a smooth
  projective variety. {I}.
\newblock {\em Inst. Hautes \'{E}tudes Sci. Publ. Math.}, 79(1):47--129, 1994.

\bibitem{waelder08}
Robert Waelder.
\newblock Equivariant elliptic genera and local {M}c{K}ay correspondences.
\newblock {\em Asian J. Math.}, 12(2):251--284, 2008.

\end{thebibliography}
\bibliographystyle{plain}
%\end{comment}
\end{document}